\documentclass[a4paper,10pt,draft]{amsart}
\usepackage{amsmath,amsthm,amssymb}
\usepackage{amsfonts}
\usepackage{latexsym}
\usepackage[dvips]{graphics}
\usepackage{amsrefs}
\usepackage{eucal}
\usepackage{mathrsfs}
\usepackage{ulem}
\usepackage{color}
\theoremstyle{plain}
\newtheorem{thm}{Theorem}[section]
\newtheorem{prop}[thm]{Proposition}
\newtheorem{lemma}[thm]{Lemma}
\newtheorem{cor}[thm]{Corollary}

\newtheorem{condition}[thm]{Condition}
 
\theoremstyle{definition}
\newtheorem{defn}[thm]{Definition}

\newtheorem{construction}[thm]{Construction}

\newtheorem{rem}[thm]{Remark}
\newtheorem{warning}[thm]{Warning}

\numberwithin{equation}{section}

 \def\subrel#1#2{\mathrel{\mathop{#2}\limits_{#1}}}

\title{On graded $\E$-rings and projective schemes 
in spectral algebraic geometry} 

\date{July 3, 2020 (version~1.1)}

\author{Mariko Ohara}
\author{Takeshi Torii}
\thanks{The second author was partially supported by
JSPS KAKENHI Grant Number JP17K05253.}
\address{Department of Mathematics, 
Shinshu University,
Matsumoto 390--8621, Japan}  
\email{primarydecomposition@gmail.com}

\address{Department of Mathematics, 
Okayama University,
Okayama 700--8530, Japan}
\email{torii@math.okayama-u.ac.jp}

\newcommand{\N}{\mathbb{N}}

\newcommand{\Z}{{\mathbb{Z}}}
\newcommand{\Spec}{\mathrm{Spec}\,}
\newcommand{\Proj}{\mathrm{Proj}\,}
\newcommand{\Hom}{\mathrm{Hom}}

\newcommand{\Cat}[1]{\mathrm{Cat}_{#1}}
\newcommand{\Ab}{\mathrm{Ab}}

\newcommand{\aperf}{\mathrm{aperf}}
\newcommand{\afg}{\mathrm{afg}}
\newcommand{\ba}{\mathrm{ba}}
\newcommand{\fg}{\mathrm{fg}}
\newcommand{\ator}{\mathrm{ator}}
\newcommand{\tor}{\mathrm{tor}}
\newcommand{\Fin}{\mathcal{F}\mathrm{in}_*}
\newcommand{\qafg}{\mathrm{qafg}}
\newcommand{\Map}{\mathrm{Map}}

\newcommand{\Mod}{\mathrm{Mod}}
\newcommand{\CAlg}{\mathrm{CAlg}}
\newcommand{\CAlgzar}{\mathrm{CAlg}^{\mathrm{Zar}}}
\newcommand{\CRing}{\mathrm{CAlg}^{\heartsuit}}
\newcommand{\CRingloc}{\mathrm{CAlg}^{\heartsuit\mathrm{Zar}}}
\newcommand{\Nil}{\mathrm{Nil}}
\newcommand{\Fun}{\mathrm{Fun}}

\newcommand{\spasce}{\operatorname{Sp}}

\newcommand{\QCoh}{\operatorname{QCoh}}

\newcommand{\Coh}{\operatorname{Coh}}
\newcommand{\colim}{\operatorname{colim}}

\newcommand{\Sp}{\mathrm{Sp}}
\newcommand{\SpDM}{\mathrm{SpDM}}

\newcommand{\Spet}{\mbox{Sp\'{e}t}}

\newcommand{\perf}{\mathrm{perf}}

\newcommand{\E}{\mathbb{E}_{\infty}}

\usepackage{enumerate}
\usepackage{array}
\usepackage[all]{xy} 
\usepackage{delarray}
\ifx\undefined\bysame
\newcommand{\bysame}{\leavemode\hbox to3em{\hrulefill}\,}
\fi
\begin{document}
\thispagestyle{empty}

\begin{abstract}
We introduce graded $\E$-rings and graded modules over them,
and study their properties.
We construct projective schemes associated
to connective $\N$-graded $\E$-rings 
in spectral algebraic geometry.
Under some finiteness conditions,
we show that the $\infty$-category of almost perfect quasi-coherent sheaves
over a spectral projective scheme $\Proj(A)$ 
associated to a connective $\N$-graded $\E$-ring $A$
can be described in terms of
$\Z$-graded $A$-modules.
\end{abstract}

\maketitle



\section{Introduction}
In stable homotopy theory, 
$\E$-rings are regarded
as counterparts of commutative rings. 
In \cite{EKMM}
Elmendorf, Kriz, Mandell, and May 
have developed the theory of rings and modules in spectra.
Lurie has formulated the theory of rings and modules
in spectra in the setting of $\infty$-categories
in \cite{HA}.
A generalization of algebraic geometry
in the setting of simplicial model categories
has been formulated and studied 
by To\"{en} and Vezzosi in \cites{TV1,TV2}.
In \cites{DAG5,DAG7,DAG8,SAG} 
Lurie has given another formulation of a generalization of 
algebraic geometry based on $\E$-rings,
which is called spectral algebraic geometry.

The aim of this paper is to introduce projective schemes
associated to connective $\N$-graded $\E$-rings,
where $\N$ is the set of nonnegative integers,
and to prove an analogue of ``the Serre theorem'' 
in spectral algebraic geometry. 

We recall ``the Serre theorem''
in classical algebraic geometry,
which describes the category of coherent sheaves on a projective scheme. 
Let $R$ be an $\N$-graded commutative ring.
We assume that $R_0$ is Noetherian
and $R$ is generated by finitely many elements in $R_1$
as an $\N$-graded commutative ring over $R_0$. 
We consider a projective scheme $\Proj(R)$ and
the category of coherent sheaves on it.
There is a functor 
$\widetilde{(-)}: \Mod_R^{\heartsuit}(\Z)
\to \QCoh(\Proj(R))$,
where $\Mod_R^{\heartsuit}(\Z)$ is the category of $\Z$-graded
$R$-modules and $\QCoh(\Proj(R))$ is the category
of quasi-coherent sheaves on $\Proj(R)$. 
If a $\Z$-graded $R$-module $M$ is
finitely generated,
then $\widetilde{M}$ is a coherent sheaf on $\Proj(R)$. 
The Serre theorem says that 
the functor $\widetilde{(-)}$ induces an
equivalence between
the category of coherent sheaves on $\Proj(R)$ 
and the quotient category of the
abelian category of finitely generated 
$\Z$-graded $R$-modules by 
the Serre subcategory of torsion $R$-modules
(cf.~\cite{Stack}*{Proposition~30.14.4}).
We say that it is the Serre theorem 
since the equivalence
for projective spaces
goes back to Serre~\cite{Serre}. 

In order to generalize the Serre theorem to spectral algebraic geometry, 
we develop the theory of
graded $\E$-rings and graded modules over them. 
A notion of graded spectra in the setting of
$\infty$-categories
have appeared in \cite{K2}*{\S2.4}.
If we regard $\Z$ as a symmetric monoidal category 
in which the objects are integers,
the morphisms are identities,
and the tensor product is given by the addition,
then we obtain a symmetric monoidal structure
on the $\infty$-category of $\Z$-graded spectra
by using Day convolution in the setting
of $\infty$-categories developed by \cite{Saul}.
Based on this symmetric monoidal structure,
we construct an $\infty$-category of
$\N$-graded $\E$-rings and
an $\infty$-category of
$\Z$-graded modules over an $\N$-graded $\E$-ring. 

Now, we recall the definition of spectral schemes in \cite{SAG}.
A spectrally ringed space is a pair $(X,\mathcal{O})$
of a topological space $X$ and 
a sheaf of $\E$-rings on $X$.
As in classical algebraic geometry,
a spectral scheme 
is defined to be a spectrally ringed space
$(X,\mathcal{O})$
which is locally equivalent to an affine
spectral scheme. 
If $(X,\mathcal{O})$ is a spectral scheme,
then the pair $(X,\pi_0\mathcal{O})$
is an ordinary scheme,
which is called the underlying scheme of $(X,\mathcal{O})$.

In this paper we construct spectral projective
schemes associated to 
connective $\N$-graded $\E$-rings.
For a connective $\N$-graded $\E$-ring $A$,
we have an ordinary projective scheme $\Proj(\pi_0(A))$
associated to the $\N$-graded commutative ring $\pi_0(A)$.
We construct a sheaf of $\E$-rings
on the underlying topological space 
of $\Proj(\pi_0(A))$ by assigning
the degree $0$ part $\E$-ring $A[f^{-1}]_0$
of the localization $A[f^{-1}]$
to an elementary open subset $D_+(f)$
for homogeneous $f\in\pi_0(A)$ of positive degree.
In this way we construct a spectral projective scheme
$\Proj(A)$ whose underlying scheme is $\Proj(\pi_0(A))$.

For a spectral scheme $\mathsf{X}$,
there is an $\infty$-category $\QCoh(\mathsf{X})$ of
quasi-coherent sheaves on it.
In the case of a spectral projective scheme $\Proj(A)$,
we construct 
a symmetric monoidal functor
\[ \widetilde{(-)}: \Mod_A(\Z)\longrightarrow \QCoh(\Proj(A)),\]
where $\Mod_A(\Z)$ is the $\infty$-category
of $\Z$-graded $A$-modules.

A quasi-coherent sheaf $\mathscr{F}$ on a spectral scheme $\mathsf{X}$ is 
said to be almost perfect if 
it corresponds to an almost perfect $B$-module
when it is restricted to each affine open subscheme 
$\Spec(B)$ of $\mathsf{X}$.
We think that the notion of almost perfect quasi-coherent
sheaves is a generalization in spectral algebraic geometry
for that of coherent sheaves on ordinary schemes.

One of our goals in this paper is to describe 
the $\infty$-category $\QCoh(\Proj A)^{\aperf}$ 
of almost perfect quasi-coherent sheaves
on $\Proj(A)$ in terms of $\Z$-graded $A$-modules
under some finiteness conditions.
We introduce an $\infty$-category
$\Mod_A^{\afg}(\Z)$ of 
almost finitely generated $\Z$-graded $A$-modules.
When the functor $\widetilde{(-)}$
is restricted to $\Mod_A^{\afg}(\Z)$,
it factors through $\QCoh(\Proj(A))^{\aperf}$.
Let $\Mod_A^{\ator}(\Z)$ be
the full subcategory of $\Mod_A^{\afg}(\Z)$
spanned by those objects $M$
such that $\widetilde{M}\simeq 0$.
We give a characterization of $M\in \Mod_A^{\ator}(\Z)$
in terms of the $\Z$-graded $\pi_0(A)$-modules
$\pi_n(M)$ for $n\in\Z$.

The following is our main result.  

\begin{thm}[Theorem~\ref{n3}]
Let $A$ be a connective Noetherian
$\N$-graded $\E$-ring satisfying 
Condition~\ref{condition:N-graded-E-infinity-ring}. 
Then, the functor $\widetilde{(-)}: \Mod_A(\Z) \to 
\QCoh(\Proj(A))$
induces an equivalence
\[ \Mod_A^{\afg}(\Z)/  \Mod_A^{\ator}(\Z)
   \stackrel{\simeq}{\longrightarrow}
   \QCoh(\Proj(A))^{\aperf}, \]
of small symmetric monoidal stable $\infty$-categories,
where the left hand side is the Verdier quotient
of small stable $\infty$-categories.
\end{thm}

The organization of this paper is as follows:
In \S\ref{section:rings-module-graded-spectra}
we review $\infty$-categories of graded spectra and
study commutative rings and modules in them.
In \S\ref{section:spectral-projective-scheme}
we introduce spectral projective schemes
associated to connective $\N$-graded $\E$-rings
and study their properties.
In \S\ref{section:quasi-coheren-sheaf}
we study the $\infty$-category
of quasi-coherent sheaves on projective spectral schemes.
We construct the functor $\widetilde{(-)}$
and show that it is a localization functor.
In \S\ref{section:almost_perfect_projective_scheme}
we give a characterization of
modules of global sections for
almost perfect quasi-coherent sheaves
on spectral projective schemes.
Using this characterization,
we prove the Serre theorem in spectral algebraic geometry.

\begin{warning}
In this paper we think ordinary categories
are a special kind of $\infty$-categories.
Thus, we does not distinguish notationally
between ordinary categories and their nerves.
\end{warning}

\subsection*{Acknowledgements}
The first author would like to express 
her thanks to Professor David Gepner for his
valuable advice on grading. 
She also would
like to express her thanks to Professor Gon\c{c}alo Tabuada 
for a lot of advice. 
She had an opportunity to talk with
them during a program 
``$K$-theory and related fields'' 
at Hausdorff Institute in Bonn. 

\section{Commutative rings and modules in graded spectra}
\label{section:rings-module-graded-spectra}

In this section we study commutative rings and modules in
$\infty$-categories of graded spectra.
In \S\ref{subsection:graded-spectra}
we review $\infty$-categories of graded spectra.
We recall that an $\infty$-category of
spectra graded by a commutative monoid
admits a symmetric monoidal structure.
In \S\ref{subsection:commutative-ring-module-graded-spectra}
we construct $\infty$-categories
of graded $\E$-rings and 
those of graded modules over graded $\E$-rings.
In \S\ref{subsection:localization-graded-ring}
we construct localizations of graded $\E$-rings
and study their mapping spaces.

\subsection{Graded spectra}
\label{subsection:graded-spectra}

In this subsection we review
$\infty$-categories of 
graded spectra
discussed in \cite{K2}*{\S2.4}
and study their properties.

First, we introduce 
a notation of indexed objects
in an $\infty$-category.

\begin{defn}
Let $S$ be a set.
We regard $S$ as a category
in which the objects are elements of $S$
and the morphisms are identity maps.
Let $\mathcal{C}$ be an $\infty$-category.
We define 
an $S$-indexed object
in $\mathcal{C}$
to be a functor $S\to \mathcal{C}$.
A morphism between 
$S$-indexed objects
is a natural transformation of functors.
We denote by 
\[ \mathcal{C}(S) \]
the $\infty$-category
$\Fun(S,\mathcal{C})$
of functors from $S$ to $\mathcal{C}$.
We call $\mathcal{C}(S)$
the $\infty$-category of 
$S$-indexed objects
in $\mathcal{C}$.
For 
an $S$-indexed object
$X$ and $s\in S$,
we set $X_s=X(s)$ for simplicity.
\end{defn}

Let $\Sp$ be the $\infty$-category of spectra.
We consider the $\infty$-category $\Sp(S)$
of $S$-indexed spectra.

\begin{prop}\label{prop:graded-spectra-stable-presentable}
The $\infty$-category $\Sp(S)$
is stable and presentable.
\end{prop}

\begin{proof}
Since $\Sp$ is a presentable stable $\infty$-category,
$\Sp(S)$ is stable 
by \cite{HA}*{Proposition~1.1.3.1},
and presentable
by \cite{HT}*{Proposition~5.5.3.6}.
\end{proof}

Next, we consider a $t$-structure
on $\Sp(S)$. 
We have a $t$-structure on $\Sp$
by \cite{HA}*{Proposition~1.4.3.6}.
Let $(\Sp{}_{\ge 0},\Sp{}_{\le 0})$
be the $t$-structure on $\Sp$.
The full subcategories $\Sp{}_{\ge 0}$
and $\Sp{}_{\le 0}$ are
spanned by spectra $X$
such that $\pi_n(X)\simeq 0$ for $n<0$ and $n>0$,
respectively. 
The $t$-structure on $\Sp$
is both left complete and right complete,
and the heart $\Sp^{\heartsuit}$
is canonically equivalent to 
$\Ab$,
where $\Ab$ is the category of abelian groups.

We define a pair of full subcategories 
\[ (\Sp(S)_{\ge 0},\Sp(S)_{\le 0}) \]
of $\Sp(S)$ to be
$(\Fun(S,\Sp{}_{\ge 0}),\Fun(S,\Sp{}_{\le 0}))$.
By \cite{HA}*{Proposition~1.4.3.6},
we easily obtain the following proposition. 

\begin{prop}\label{prop:t-structure-graded-spectra}
The pair $(\Sp(S)_{\ge 0},\Sp(S)_{\le 0})$
of full subcategories 
determines an accessible $t$-structure
on $\Sp(S)$.
The $t$-structure on $\Sp(S)$
is both left complete and right complete,
and the heart $\Sp(S)^{\heartsuit}$
is canonically equivalent to 
$\Ab(S)$.
\end{prop}

\begin{defn}
We say that 
an $S$-indexed spectrum 
is connective if
it belongs to $\Sp(S)_{\ge 0}$.
\end{defn}

For any integer $n$, 
the functor $\pi_n: \Sp\to\Ab$ induces a functor 
\[ \pi_n: \Sp(S)\longrightarrow \Ab(S).\]
An $S$-indexed spectrum $X$ 
is connective if and only if
$\pi_n(X_s)\simeq 0$ for all $s\in S$ and $n<0$.
If $a\in \pi_n(X_s)$,
then we say that $a$ is a homogeneous element of $\pi_n(X)$
of degree $s$.

Let $T$ be a subset of $S$.
We regard $T$ as a full subcategory of $S$.
The inclusion map $i: T\hookrightarrow S$ induces
a restriction functor
\[ i^*: \Sp(S)\longrightarrow \Sp(T).\]
By a left Kan extension along $i$,
we obtain a left adjoint
\[ i_!: \Sp(T)\longrightarrow \Sp(S) \]
to $i^*$.
Note that $i_!$ is also a right adjoint to $i^*$.
Since the unit ${\rm id} \to i^*\circ i_!$ of the adjunction
is an equivalence,
$i_!$ is fully faithful.
In the following of this paper,
we identify $\Sp(T)$ with 
the full subcategory of $\Sp(S)$
spanned by objects in the essential image of $i_!$.

Now, we suppose that $G=(G,+,0)$ is a commutative monoid.
We consider a symmetric monoidal structure
on the $\infty$-category $\Sp(G)$ of 
$G$-indexed spectra.
See \cite{HA}*{\S2.1}
and \cite{Groth}*{\S4.4} 
for symmetric monoidal $\infty$-categories.
In particular,
see \cite{HA}*{Definition~2.0.0.7}
and \cite{Groth}*{Definition~4.27}
for the definition of symmetric monoidal $\infty$-categories.

We regard $G$ as a strict symmetric monoidal category
in which the objects are elements of $G$,
the morphisms are identity maps, and
the tensor product is given by the addition $+$ of $G$.
Since $\Sp$ is a symmetric monoidal $\infty$-category
in which the tensor product preserves
small colimits separately in each variable,
we obtain the following proposition
by the theory of Day convolution in the setting of $\infty$-categories
developed by \cite{Saul}.
For another construction
of symmetric monoidal structure
on ${\rm Sp}(G)$, see \cite{K2}*{\S2.3}.

\begin{prop}[{cf.~\cite{Saul}*{Proposition~2.11
and Lemma~2.13}, and
\cite{K2}*{Corollary~2.3.9 and Remark~2.3.10}}]
\label{prop:symmetric-monoidal-structure-graded-spectra}
The $\infty$-category $\Sp(G)$ of 
$G$-indexed spectra
admits a symmetric monoidal structure
by Day convolution,
in which the tensor product preserves 
small colimits separately in each variable.
\end{prop}

In the following of this paper we regard $\Sp(G)$
as an $\infty$-category equipped with
the symmetric monoidal structure given by
Proposition~\ref{prop:symmetric-monoidal-structure-graded-spectra}.
We call an object of $\Sp(G)$
a $G$-graded spectrum.

For $G$-graded spectra $X$ and $Y$,
we denote by $X\otimes Y$
the tensor product of $X$ and $Y$
in $\Sp(G)$. 
Note that 
\[ (X\otimes Y)_g\simeq
   \bigoplus_{g'+g''=g}X_{g'}\otimes Y_{g''}\]
for $g\in G$ and that 
the unit object $I$
is given by
\[ I_g\simeq\left\{
    \begin{array}{ll}
      \mathbb{S}   & \mbox{\rm if $g=0$},\\
      0   & \mbox{\rm if $g\neq 0$},\\
    \end{array}\right. \]
where $\mathbb{S}$ is the sphere spectrum.

\begin{defn}[{cf.~\cite{K2}*{Example~2.4.8}}]
\label{definition:twisting-spectra}
Suppose that $G$ is a commutative group.
For $X\in\Sp(G)$ and $g\in G$,
we define a twisting
\[ X(g)\in \Sp(G) \]
by $X(g)_{g'}\simeq X_{g+g'}$ for $g'\in G$.
We define a twisting functor
\[ (g): \Sp(G)\to \Sp(G) \]
by assigning to $X\in\Sp(G)$
the twisting $X(g)\in \Sp(G)$.
Since $\Sp(G)$ is a symmetric monoidal $\infty$-category,
we can regard $\Sp(G)$ as tensored over itself.
In this case the twisting functor
$(g): \Sp(G)\to \Sp(G)$ is $\Sp(G)$-linear
(see \cite{HA}*{Definition~4.6.2.7}
for the definition of $\mathcal{C}$-linear functors).
\end{defn}

On the $\infty$-category $\Sp(G)$
of $G$-graded spectra,
we have the $t$-structure 
by Proposition~\ref{prop:t-structure-graded-spectra}
and 
the symmetric monoidal structure
by Proposition~\ref{prop:symmetric-monoidal-structure-graded-spectra}.
We have the following proposition
on compatibility of the $t$-structure and
the symmetric monoidal structure on $\Sp(G)$.

\begin{prop}\label{prop:compatibility-t-structure-symmetric-monoidal}
The full subcategory
$\Sp(G)_{\ge 0}$ contains the unit object and
is closed under tensor products.
\end{prop}

Let $H$ be a submonoid of $G$.
We regard $H$ as a full subcategory of $G$
and let $i: H\hookrightarrow G$ be
the inclusion functor.
It is easy to see that 
the fully faithful functor
$i_!: \Sp(H)\longrightarrow \Sp(G)$
is symmetric monoidal
(see \cite{HA}*{Definition~2.1.3.7} and
\cite{Groth}*{Definition~4.31}
for the definition of symmetric monoidal functors).
In other words,
the full subcategory $\Sp(H)$
contains the unit object and is closed
under tensor products
of $\Sp(G)$.
On the other hand,
the restriction functor 
$i^*: \Sp(G)\to\Sp(H)$ 
is lax symmetric monoidal in general
(see \cite{HA}*{Definition~2.1.2.7}
and \cite{Groth}*{Definition~4.31}
for the definition of lax symmetric monoidal functors).

\begin{defn}
We define a functor
\[ (-)_0: \Sp(G)\longrightarrow \Sp \] 
to be the restriction functor $j_0^*$
induced by the inclusion
$j_0: \{0\}\hookrightarrow G$,
where we identify $\Sp(\{0\})$
with $\Sp$.
\end{defn}

\subsection{Commutative rings and modules in graded spectra}
\label{subsection:commutative-ring-module-graded-spectra}

In this subsection we will define
commutative rings and modules
in $\infty$-categories of graded spectra
and study their properties.

Let $G$ be a commutative monoid.
We consider the $\infty$-category $\Sp(G)$
of $G$-graded spectra.
Recall that $\Sp(G)$ is a 
symmetric monoidal $\infty$-category 
by Proposition~\ref{prop:symmetric-monoidal-structure-graded-spectra}.

First, we consider commutative rings 
in the $\infty$-category of $G$-graded spectra.
We call a commutative algebra object 
in $\Sp(G)$ a $G$-graded $\E$-ring.
We define
\[  \CAlg(G) \]
to be the $\infty$-category
$\CAlg(\Sp(G))$
of $G$-graded $\E$-rings.
We say that a $G$-graded $\E$-ring is connective 
if its underlying $G$-graded spectrum is connective. 

\begin{prop}
The $\infty$-category
$\CAlg(G)$ is presentable.
\end{prop}

\begin{proof}
The proposition follows from
\cite{HA}*{Corollary~3.2.3.5}.
\end{proof}

Let $H$ be a submonoid of $G$
and let $i: H\hookrightarrow G$ be the inclusion map.
Recall that we have an adjunction
\[ i_!: \Sp(H)\rightleftarrows \Sp(G): i^* \]
of $\infty$-categories,
where the left adjoint $i_!$ is fully faithful
and symmetric monoidal,
and the right adjoint is lax symmetric monoidal.   
By \cite{HA}*{Remark~7.3.2.13},
this induces an adjunction
\[ i_!: \CAlg(H)\rightleftarrows \CAlg(G):i^*\]
of $\infty$-categories,
where the left adjoint $i_!$ is fully faithful.
We understand that an $H$-graded $\E$-ring is a
$G$-graded $\E$-ring $A$ satisfying $A_g\simeq 0$
for $g\not\in H$. 

\begin{defn}\label{definition:restriction-CAlg}
We define a functor
\[ (-)_0: \CAlg(G)\longrightarrow
          \CAlg \]
to be the restriction functor $j_0^*$  
induced by the inclusion $j_0: \{0\}\hookrightarrow G$,
where we identify $\CAlg(\{0\})$
with the $\infty$-category $\CAlg$ of $\E$-rings.
\end{defn}

Next, we consider graded modules 
over graded $\E$-rings. 
Let $A$ be a $G$-graded $\E$-ring.
We call an $A$-module object in $\CAlg(G)$ 
a $G$-graded $A$-module.
We define
\[ \Mod_A(G) \]
to be the $\infty$-category
$\Mod_A(\Sp(G))$
of $G$-graded $A$-modules.

\begin{prop}\label{prop:module-category-stable-presentable}
The $\infty$-category $\Mod_A(G)$ is stable
and presentable.
\end{prop}

\begin{proof}
By \cite{HA}*{Corollary~4.5.1.6},
we have an equivalence of $\infty$-categories
between $\Mod_A(G)$ and 
the $\infty$-category $\mathrm{LMod}_A(\Sp(G))$
of left $A$-module objects in $\Sp(G)$.
By Proposition~\ref{prop:symmetric-monoidal-structure-graded-spectra},
$\Sp(G)$ is a symmetric monoidal $\infty$-category
in which the tensor product preserves
small colimits separately in each variable.
Since $\Sp(G)$ is stable and presentable
by Proposition~\ref{prop:graded-spectra-stable-presentable},
$\mathrm{LMod}_A(\Sp(G))$ is stable
by \cite{HA}*{Proposition~7.1.1.4},
and presentable
by \cite{HA}*{Corollary~4.2.3.7(1)}.
\end{proof}

\begin{defn}\label{definition:twisting-module-spectra}
Suppose that $G$ is a commutative group.
Recall that we have defined the twisting functor
$(g): \Sp(G)\to\Sp(G)$ for $g\in G$
in Definition~\ref{definition:twisting-spectra}.
Since it is $\Sp(G)$-linear,
the twisting $M(g)$ supports a natural $A$-module
structure for $M\in\Mod_A(G)$.
We define a twisting functor
\[ (g): \Mod_A(G)\longrightarrow \Mod_A(G) \]
by assigning to $M\in\Mod_A(G)$
the twisting $M(g)\in\Mod_A(G)$.
\end{defn}

Now we shall show that 
$\Mod_A(G)$
admits a symmetric monoidal structure.

\begin{prop}
The $\infty$-category
$\Mod_A(G)$ admits
a symmetric monoidal structure
in which the tensor product 
is given by the relative tensor product $\otimes_A$
and the unit object is $A$.
The relative tensor product $\otimes_A$
preserves small colimits separately
in each variable.
\end{prop}

\begin{proof}
By Propositions~\ref{prop:graded-spectra-stable-presentable}
and \ref{prop:symmetric-monoidal-structure-graded-spectra},
$\Sp(G)$
is a presentable symmetric monoidal $\infty$-category
in which 
the tensor product preserves 
small colimits separately in each variable.
The proposition follows from
\cite{HA}*{Theorem~4.5.2.1 and Corollary~4.4.2.15}.
\end{proof}

We define an $\infty$-category
\[ \CAlg_A(G) \]
to be $\CAlg(\Mod_A(G))$.
We call $\CAlg_A(G)$
the $\infty$-category of $G$-graded $\E$-rings over $A$.
Note that there is an equivalence
\[ \CAlg_A(G)\simeq \CAlg(G)_{A/} \]
of $\infty$-categories
by \cite{HA}*{Corollary~3.4.1.7}.

Let $H$ be a submonoid of $G$.
The restriction functor $i^*: \Sp(G)\to \Sp(H)$
induces a restriction functor
\[ i^*: \Mod_A(G)\longrightarrow \Mod_B(H),\]
where $B=i^*A\in\CAlg(H)$. 
By the adjunction
$i_!:\CAlg(H)\rightleftarrows \CAlg(G):i^*$,
there is a map $B\to A$ of $G$-graded $\E$-rings,
where we identify $B$ with $i_!B$.
Then we have an adjunction
\[ i_!: \Mod_B(H)\rightleftarrows
   \Mod_A(G): i^* \]  
of $\infty$-categories,
where the left adjoint $i_!$ 
is the scalar extension functor $A\otimes_B(-)$.
Note that the left adjoint $i_!$
is symmetric monoidal and 
the right adjoint $i^*$ is lax symmetric monoidal.

\begin{defn}
We define a functor
\[ (-)_0: \Mod_A(G)\longrightarrow \Mod_{A_0} \]
to be the restriction functor $j_0^*$
induced by the inclusion 
$j_0: \{0\}\hookrightarrow G$,
where $\Mod_{A_0}$ is the $\infty$-category
of module spectra over the $\E$-ring $A_0$. 
\end{defn}

Now, we shall show that 
the adjunction 
$i_!: \Mod_B(H)\to\Mod_A(G): i^*$
gives an equivalence
of $\infty$-categories under some conditions.
We suppose that $G\cong H\oplus \Z$,
where $\Z$ is the set of integers. 
In this case, for any $H$-graded $B$-module $N$,
the unit map $N\to i^*(A\otimes_B N)$
of the adjunction is an equivalence.
Hence the left adjoint
$i_!$ is fully faithful.

\begin{prop}\label{prop:periodic-equivalence}
Let $G$ be a commutative monoid and 
let $H$ be a submonoid of $G$ 
such that $G\cong H\oplus\Z$.
Let $A$ be a $G$-graded $\E$-ring.
We assume that there exists $a\in \pi_0(A_g)$
such that $a$ is invertible in the $G$-graded commutative 
ring $\pi_0(A)$,
where $g\in G$ corresponds to
$(0,1)$ under the isomorphism $G\cong H\oplus\Z$.
Then the restriction functor induces an equivalence
\[ i^*: \Mod_A(G)\stackrel{\simeq}{\longrightarrow} \Mod_{B}(H) \]
of symmetric monoidal stable $\infty$-categories.
\end{prop}

\begin{proof}
Since we have the adjunction
$i_!: \Mod_B(H)\rightleftarrows \Mod_{A}(G):i^*$
and the left adjoint $i_!$ is fully faithful,
it suffices to show that the counit
of the adjunction is an equivalence.
Note that, in order to show that
a map $X\to Y$ between $G$-graded $A$-modules
is an equivalence, it suffices
to show that the restrictions
$i^*(X(ng))\to i^*(Y(ng))$ are equivalences
for all $n\in\Z$.

For any $G$-graded $A$-module $M$,
we set $N=i^*M$.
There is a commutative diagram
\[ \begin{array}{ccc}
    A\otimes_{B}N & \longrightarrow & M\\[1mm]
    \mbox{$\scriptstyle a^n\otimes {\rm id}$}
    \bigg\downarrow
    \phantom{\mbox{$\scriptstyle a^n\otimes {\rm id}$}}
     & & 
    \phantom{\mbox{$\scriptstyle a^n$}}
    \bigg\downarrow
    \mbox{$\scriptstyle a^n$}\\[3mm]
    A(ng)\otimes_{B}N & \longrightarrow &
    M(ng),\\
   \end{array}\]
where the vertical arrows are equivalences
by the assumption.
Applying the functor $i^*$
to the above diagram,
we obtain a commutative diagram
\[ \begin{array}{ccc}
    B\otimes_{B}N & \longrightarrow & N\\[1mm]
    \mbox{$\scriptstyle a^n\otimes {\rm id}$}
    \bigg\downarrow
    \phantom{\mbox{$\scriptstyle a^n\otimes {\rm id}$}}
     & & 
    \phantom{\mbox{$\scriptstyle a^n$}}
    \bigg\downarrow
    \mbox{$\scriptstyle a^n$}\\[3mm]
    i^*(A(ng)\otimes_{B}N) & \longrightarrow &
    i^*(M(ng)),\\
   \end{array}\]
where the vertical arrows are equivalences.
Since the top horizontal arrow is an equivalence,
the bottom horizontal arrow is also an equivalence.
This shows that
the map
$i^*((A\otimes_{B}N)(ng))\to i^*(M(ng))$ is an equivalence
for all $n\in\Z$.
Thus, 
the counit map $i_!i^*M\simeq A\otimes_BN\to M$
is an equivalence. 
\end{proof}

\begin{cor}\label{cor:periodic-equivalence}
Let $A$ be a $\Z$-graded $\E$-ring.
Suppose that there exists $a\in \pi_0(A_1)$ such that
$a$ is invertible in the $\Z$-graded commutative ring $\pi_0(A)$.
Then the restriction functor
induces an equivalence
\[ i^*: \Mod_A(\Z)\stackrel{\simeq}{\longrightarrow} 
        \Mod_{A_0} \]
of symmetric monoidal stable $\infty$-categories.
\end{cor}

\begin{rem}
By \cite{K2}*{p.17, the second
clause before Definition~2.4.2},
the homomorphism $G\to \{0\}$
induces a symmetric monoidal functor
\[ \mathrm{Und}: \Sp(G)\longrightarrow \Sp.\]
For $X\in\Sp(G)$,
we call $\mathrm{Und}(X)$ the underlying spectrum of $X$,
which is given by
\[ \mathrm{Und}(X)\simeq\bigoplus_{g\in G}X_g. \]
Since $\mathrm{Und}$ is symmetric monoidal,
it induces a symmetric monoidal functor
\[ \mathrm{Und}: \Mod_A(G)\longrightarrow
                 \Mod_{\mathrm{Und}(A)}.\]
\end{rem}

\subsection{Localizations of graded $\E$-rings}
\label{subsection:localization-graded-ring}

In this subsection 
we construct localizations of $G$-graded $\E$-rings.
For a $G$-graded $\E$-ring $A$,
we study mapping spaces from localizations
of $A$ to any $G$-graded $\E$-rings over $A$.
Throughout this subsection
we assume that $G$ is a commutative group.

Let $A$ be a $G$-graded $\E$-ring
and let $a\in \pi_0(A)$ be homogeneous of degree $g\in G$.
We regard $a$ as a morphism
$a: A\to A(g)$ of $G$-graded $A$-modules.
Since $\Mod_A(G)$ is a presentable $\infty$-category
by Proposition~\ref{prop:module-category-stable-presentable},
there exists a localization functor
\[ L: \Mod_A(G)\longrightarrow \Mod_A(G) \]
with respect to the map $a: A\to A(g)$.
By definition,
$M\in\Mod_A(G)$ is $L$-local
if and only if the multiplication map
$a: M\to M(g)$ is an equivalence.
As in the nongraded case,
$L(M)$ is equivalent to $M[a^{-1}]$,
where $M[a^{-1}]$ is
a colimit of the sequence
\[ M\stackrel{a}{\to}M(g)\stackrel{a}{\to}M(2g)
   \stackrel{a}{\to}\cdots \]
in $\Mod_A(G)$.
In particular,
$L$ is a smashing localization given by
$L(-)\simeq A[a^{-1}]\otimes_A(-)$.

The localization map $l:A\to A[a^{-1}]$
is an idempotent object of $\Mod_A(G)$
in the sense of \cite{HA}*{Definition~4.8.2.1}.
By \cite{HA}*{Proposition~4.8.2.7},
$L$ is compatible with the symmetric monoidal
structure on $\Mod_A(G)$
(see also \cite{HA}*{Proposition~2.2.1.9}).
Hence,
$L$ is a lax symmetric monoidal functor.
We can regard 
$L(A)\simeq A[a^{-1}]$ as a commutative monoid object
of $\Mod_A(G)$,
and $l: A\to A[a^{-1}]$ a morphism in $\CAlg(G)$.
We obtain an adjunction
\[ l_!: \Mod_A(G)\rightleftarrows \Mod_{A[a^{-1}]}(G): l^*,\]  
where the left adjoint $l_!$ is 
a symmetric monoidal
functor given by $M\mapsto A[a^{-1}]\otimes_A M$,
and the right adjoint $l^*$ is 
a fully faithful lax symmetric monoidal functor. 

By \cite{HA}*{Remark~7.3.2.13},
this adjunction induces an adjunction
\[ l_!: \CAlg_A\rightleftarrows \CAlg_{A[a^{-1}]}: l^*, \]
where the right adjoint $l^*$ is fully faithful.
Hence $l_!: \CAlg_A(G)\to\CAlg_{A[a^1]}(G)$
is a localization functor.

\begin{rem}
We notice that the following conditions for $B\in\CAlg_A(G)$
are equivalent:
\begin{enumerate}
\item
The $G$-graded $\E$-ring $B$ is $l_!$-local.
\item 
The $G$-graded $\E$-ring 
$B$ belongs  to the essential image of $l^*$.
\item
The unit map $B\to l^*l_!B$ of the adjunction
is an equivalence in $\CAlg_A(G)$.
\item
The localization map $B\to B[a^{-1}]$
is an equivalence in $\Mod_A(G)$.
\item
The element $a$ is invertible
in the $G$-graded commutative ring
$\pi_0(B)$.
\end{enumerate}
\end{rem}

\begin{lemma}\label{lemma:mapping_space_graded_rings}
Let $A$ be a $G$-graded $\E$-ring
and let $a\in\pi_0(A)$ be a homogeneous element.
For any $G$-graded $\E$-ring $B$ over $A$, 
the mapping space
\[ \Map_{\CAlg_A(G)}(A[a^{-1}],B)\]
is $(-1)$-truncated.
It is nonempty if and only if
$a$ is invertible in the $G$-graded commutative ring
$\pi_0(B)$.
\end{lemma}

\begin{proof}
Clearly the mapping space is empty
if $a$ is not invertible in $\pi_0(B)$. 
If $a$ is invertible in $\pi_0(B)$,
then $B$ is  $l_!$-local.
In this case
we have an equivalence
\[ \Map_{\CAlg_A(G)}(A[a^{-1}],B)\stackrel{\simeq}{\longrightarrow}
   \Map_{\CAlg_A(G)}(A,B) \]
since there is an equivalence
$l_!A\simeq A[a^{-1}]$ in $\CAlg_A(G)$.
The lemma follows from the fact that 
the mapping space
$\Map_{\CAlg_A(G)}(A,B)$ is contractible.
\end{proof}

\begin{cor}\label{cor:mapping_space_localizaed-graded_rings}
Let $A$ be a $G$-graded $\E$-ring
and let $a\in\pi_0(A)$ be a homogeneous element.
For any $G$-graded $\E$-ring $B$,
the induced map 
\[ \Map_{\CAlg(G)}(A[a^{-1}],B)\longrightarrow
   \Map_{\CAlg(G)}(A,B)\]
restricts to a homotopy equivalence of 
$\Map_{\CAlg(G)}(A[a^{-1}],B)$ with the summand of
$\Map_{\CAlg(G)}(A,B)$ spanned by those maps $A\to B$
which carry $a\in\pi_0(A)$ to an invertible element of
$\pi_0(B)$.
\end{cor}

\begin{proof}
Let $f: A\to B$ be a map of $G$-graded 
$\E$-rings.
The mapping space $\Map_{\CAlg_A(G)}(A[a^{-1}],B)$
is equivalent to the fiber of the map
$\Map_{\CAlg(G)}(A[a^{-1}],B)\to\Map_{\CAlg(G)}(A,B)$
at $f$
since there is an equivalence $\CAlg_A(G)\simeq \CAlg(G)_{A/}$.
Hence the corollary follows from 
Lemma~\ref{lemma:mapping_space_graded_rings}.
\end{proof}

We define an $\infty$-category
\[ \CAlgzar_A(G) \]
to be the full subcategory of
$\CAlg_A(G)$ spanned by
those objects of the form $A[a^{-1}]$
for some homogeneous element $a\in\pi_0(A)$.
The functor $\pi_0: \Sp(G)\to \Ab(G)$
induces a functor
\[ \pi_0: \CAlg_A(G)\longrightarrow \CRing_{\pi_0(A)}(G),\]
where $\CRing_{\pi_0(A)}(G)$
is the category of $G$-graded commutative rings
over $\pi_0(A)$.
We define $\CRingloc_{\pi_0(A)}(G)$
to be the full subcategory of $\CRing_{\pi_0(A)}(G)$
spanned by 
those objects
of the form $\pi_0(A)[a^{-1}]$
for some $a\in \pi_0(A)$.

\begin{thm}\label{thm:open-Zariski-immersions}
The functor $\pi_0: \CAlg_A(G)\to\CRing_{\pi_0(A)}(G)$
induces an equivalence
\[ \pi_0:    \CAlgzar_A(G)
   \stackrel{\simeq}{\longrightarrow}
          \CRingloc_{\pi_0(A)}(G)\]
of $\infty$-categories.
\end{thm}

\begin{proof}
By Lemma~\ref{lemma:mapping_space_graded_rings},
we see that the functor $\pi_0$ induces an equivalence
\[ \Map_{\CAlg_A(G)}(A[a^{-1}],A[b^{-1}])
   \stackrel{\simeq}{\longrightarrow}
   \Map_{\CRing_{\pi_0(A)}(G)}(\pi_0(A)[a^{-1}],\pi_0(A)[b^{-1}])\]
of mapping spaces
for any homogeneous elements $a,b\in\pi_0(A)$.
Hence
$\pi_0: 
\CAlgzar_A(G)\to
         \CRingloc_{\pi_0(A)}(G)$
is an equivalence of $\infty$-categories.
\end{proof}

\section{Spectral projective schemes
associated to $\mathbb{N}$-graded $\E$-rings}
\label{section:spectral-projective-scheme}

In this section 
we study spectral projective schemes.
As in classical algebraic geometry, 
a spectral scheme is a spectrally ringed space
which is locally equivalent to an affine
spectral scheme.
Taking $\pi_0$ of the structure sheaf,
we obtain the underlying ordinary scheme of 
a spectral scheme.
In \S\ref{subsection:definition-spectral-projective-scheme}
we introduce a spectral projective scheme $\Proj(A)$
associated to a connective $\mathbb{N}$-graded $\E$-ring $A$,
whose underlying scheme is the ordinary
projective scheme $\Proj(\pi_0(A))$.
In \S\ref{subsection:properties-projective-schemes}
we study properties of spectral projective schemes.
In particular, we show that a morphism of 
connective $\N$-graded $\E$-rings
induces a closed immersion of associated
spectral projective schemes
if it induces a surjection on $\pi_0$.
Finally, we show that 
Lurie's projective space is obtained as a spectral 
projective scheme
for some connective $\N$-graded $\E$-ring.

\subsection{Definition of spectral projective schemes associated
to connective $\N$-graded $\E$-rings}
\label{subsection:definition-spectral-projective-scheme}

In this subsection we introduce spectral projective schemes
associated to connective $\N$-graded $\E$-rings,
where $\N$ is the set of nonnegative integers.

First, we recall the definition of spectral schemes.
According to \cite{SAG}*{Definition~1.1.2.5},
a spectrally ringed space is a pair $(X,\mathcal{O}_X)$,
where $X$ is a topological space and $\mathcal{O}_X$ 
is a $\CAlg$-valued sheaf on $X$.

\begin{defn}[{cf.~\cite{SAG}*{Definition~1.1.2.8}}]
A spectral scheme is defined to be a spectrally ringed space
$(X,\mathcal{O}_X)$ which satisfies the following conditions:
\begin{enumerate}
\item[(1)] 
The underlying ringed space 
$(X,\pi_0\mathcal{O}_X)$
is an ordinary scheme.
\item[(2)] 
For each $n\in\Z$,
the sheaf $\pi_n\mathcal{O}_X$ is quasi-coherent
as a sheaf of $\pi_0\mathcal{O}_X$-modules on $X$.
\item[(3)] 
Let $U$ be an open subset of $X$ for which 
the scheme $(U,\pi_0\mathcal{O}_X|_U)$ is affine.
Then the canonical map $\pi_n(\mathcal{O}_X(U))\to
(\pi_n\mathcal{O}_X)(U)$ is an isomorphism
for every integer $n$.
\item[(4)] 
The sheaf $\pi_n\mathcal{O}_X$ vanishes for every $n<0$.
\end{enumerate}
When a spectrally ringed space $(X,\mathcal{O}_X)$ 
is a spectral scheme,
we say that $(X,\pi_0\mathcal{O}_X)$
is the underlying scheme of 
the spectral scheme $(X,\mathcal{O}_X)$.
\end{defn}

Basic examples of spectral schemes are supplied by
affine spectral schemes, 
which are associated to connective $\E$-rings.
Next, we recall the construction of affine spectral schemes.

\begin{construction}
Let $R$ be a connective $\E$-ring.
Note that we have a commutative ring $\pi_0(R)$.
We denote by 
\[ 
|\Spec(\pi_0(R))| 
\]
the underlying topological space of the ordinary
affine scheme $\Spec(\pi_0(R))$.

We construct a $\CAlg$-valued sheaf
on $|\Spec(\pi_0(R))|$.
Let $\mathcal{U}(R)$ be the partially ordered set of
all open subsets of $|\Spec(\pi_0(R))|$,
which we regard as a category.
For $r\in\pi_0(R)$,
we have an open subset $D(r)$
of $|\Spec(\pi_0(R))|$
given by 
\[ 
D(r)=\{\mathcal{P}\in\Spec(\pi_0(R))|\
r\not\in\mathcal{P}\}.
\]
We let $\mathcal{U}(R)_e$ be the 
full subcategory of $\mathcal{U}(R)$ spanned by
all open subsets of the form $D(r)$ for $r\in\pi_0(R)$.
Let $\CRing_{\pi_0(R)}$
be the category of commutative rings over $\pi_0(R)$.
We have a functor 
\[
\overline{\mathcal{O}}{}':
\mathcal{U}(R)_e^{\mathrm{op}}\to 
\CRing_{\pi_0(R)},
\]
which assigns to $D(r)$ the commutative ring
$\pi_0(R)[r^{-1}]$ over $\pi_0(R)$.
We define
$\CRingloc_{\pi_0(R)}$
to be the full subcategory of $\CRing_{\pi_0(R)}$
spanned by those objects of the form
$\pi_0(R)[r^{-1}]$ for some $r\in\pi_0(R)$.
Then we can regard $\overline{\mathcal{O}}{}'$
as a functor $\mathcal{U}(R)_e^{\mathrm{op}}\to
\CRingloc_{\pi_0(R)}$.

Let $\CAlgzar_R$
be the full subcategory of $\CAlg_R$
spanned by those objects of the form
$R[r^{-1}]$ for some $r\in\pi_0(R)$.
Using the fact that 
the functor $\pi_0: \CAlg_R\to \CAlg_{\pi_0(R)}$
induces an equivalence
$\pi_0: \CAlgzar_R\stackrel{\simeq}{\to} 
\CRingloc_{\pi_0(R)}$
by \cite{HA}*{Theorem~7.5.0.6},
we can lift the functor $\overline{\mathcal{O}}{}'$
to a functor 
\[ \mathcal{O}': \mathcal{U}(R)_e^{\mathrm{op}}\to\CAlgzar_R,\]
which assigns $R[r^{-1}]$ to $D(r)$.
We regard $\mathcal{O}'$
as a functor 
$\mathcal{O}': \mathcal{U}(R)_e^{\mathrm{op}}\to \CAlg$.
We define 
\[ 
\mathcal{O}: \mathcal{U}(R)^{\mathrm{op}}\to\CAlg 
\]
to be a right Kan extension of $\mathcal{O}'$
along the inclusion 
$\mathcal{U}(R)_e^{\mathrm{op}}\hookrightarrow
\mathcal{U}(R)^{\mathrm{op}}$.
In the same way as in the proof of 
\cite{SAG}*{Proposition~1.1.4.2},
we see that $\mathcal{O}$ is a $\CAlg$-valued sheaf on 
$|\Spec(\pi_0(R))|$
and that the spectrally ringed space
$(|\Spec(\pi_0(R))|,\mathcal{O})$
is a spectral scheme.
We say that  
$(|\Spec(\pi_0(R))|,\mathcal{O})$
is the affine spectral scheme associated
to $R$ and 
denote it by 
\[ \Spec(R). \]
Note that we have an equivalence
$\mathcal{O}(D(r))\simeq R[r^{-1}]$
for any $r\in\pi_0(R)$.
\end{construction}

Now we construct spectral projective schemes
associated to connective $\N$-graded $\E$-rings.

\begin{construction}\label{construction:spectral-projective-scheme}
For a connective $\N$-graded $\E$-ring $A$,
we have an $\N$-graded commutative ring $\pi_0(A)$.
We denote by 
\[ |\Proj(\pi_0(A))| \]
the underlying topological space
of the ordinary projective scheme $\Proj(\pi_0(A))$.

We construct a $\CAlg$-valued sheaf on
$|\Proj(\pi_0(A))|$.
Let $\mathcal{U}(A)$ be the partially ordered set 
of all open subsets of $|\Proj(\pi_0(A))|$,
which we regard as a category.
For a homogeneous element $f$ of $\pi_0(A)$ of positive degree,
we have an open subset 
$D_+(f)$ of $|\Proj(\pi_0(A))|$
given by
\[ D_+(f)=\{\mathcal{P}\in \Proj(\pi_0(A))
           |\ f\not\in\mathcal{P}\}. \]
We denote by $\mathcal{U}(A)_e$
the full subcategory of $\mathcal{U}(A)$
spanned by all open subsets of the form $D_+(f)$ 
for a homogeneous element
$f\in\pi_0(A)$ of positive degree.
We have a functor 
\[ \overline{\mathcal{O}}{}'_*: \mathcal{U}(A)_e^{\mathrm{op}}\to
   \CRing_{\pi_0(A)}(\Z), \]
which assigns to $D_+(f)$
the $\Z$-graded commutative ring $\pi_0(A)[f^{-1}]$
over $\pi_0(A)$.


Recall that $\CAlgzar_A(\Z)$
is the full subcategory of $\CAlg_A(\Z)$
spanned by those objects of the form
$A[f^{-1}]$ for some homogeneous element
$f\in\pi_0(A)$ of positive degree.
By Theorem~\ref{thm:open-Zariski-immersions},
$\CAlgzar_A(\Z)$
is equivalent to the ordinary category
$\CRingloc_{\pi_0(A)}(\Z)$.
Hence we can lift the functor $\overline{\mathcal{O}}{}'_*$
to a functor
\[ \mathcal{O}'_*: 
   \mathcal{U}(A)_e^{\mathrm{op}}\to\CAlgzar_{A}(\Z), \]
which assigns to $D_+(f)$
the $\Z$-graded $\E$-ring $A[f^{-1}]$ over $A$.
By 
regarding $\mathcal{O}'_*$ as a functor
to $\CAlg_A(\Z)$ and
composing with 
the restriction functor $(-)_0: \CAlg_A(\Z)\to \CAlg_{A_0}$
in Definition~\ref{definition:restriction-CAlg},
we obtain a functor 
\[ \mathcal{O}': \mathcal{U}(A)_e^{\mathrm{op}}\to\CAlg_{A_0}, \]
which assigns to $D_+(f)$ 
the $\E$-ring $A[f^{-1}]_0$
over $A_0$.
We regard $\mathcal{O}'$ as a functor
$\mathcal{O}': \mathcal{U}(A)_e^{\mathrm{op}}\to\CAlg$.
We define a functor 
\[ \mathcal{O}:
   \mathcal{U}(A)^{\mathrm op}\to \CAlg \]
to be a right Kan extension of
$\mathcal{O}'$
along the inclusion
$\mathcal{U}(A)_e^{\mathrm{op}}\to
\mathcal{U}(A)^{\mathrm{op}}$.
\end{construction}

\begin{lemma}\label{lemma:construction-structure-sheaf}
The functor $\mathcal{O}$ 
is a $\CAlg$-valued sheaf on 
$|\Proj(\pi_0(A))|$.
\end{lemma}

\begin{proof}
It suffices to show that $\mathcal{O}$
satisfies condition (2) of 
\cite{SAG}*{Proposition~1.1.4.4}.
Let $D_+(f_1),\ldots,D_+(f_n)$
be open subsets belonging to $\mathcal{U}(A)_e$
such that $D(f)=\bigcup_{1\le i\le n}D_+(f_i)$
also belongs to $\mathcal{U}(A)_e$.
Set $R=\pi_0(A)[f^{-1}]_0$.
We have a homeomorphism
$D_+(f)\cong |\Spec(R)|$.
By the construction of the functor $\mathcal{O}'$,
we see that the restriction $\mathcal{O}'|_{D_+(f)}$
can be identified with 
the functor $\mathcal{U}(R)_e^{\mathrm{op}}\to
\CAlg$ which assigns to $D(r)$
the $\E$-ring $R[r^{-1}]$.
Hence the restriction
$(D_+(f),\mathcal{O}|_{D_+(f)})$
is equivalent to the affine spectral 
scheme $\Spec(R)$.
In particular,
we see that condition (2) of 
\cite{SAG}*{Proposition~1.1.4.4}
is satisfied for $D_+(f_1),\ldots,D_+(f_n)$.
This completes the proof.
\end{proof}

By Lemma~\ref{lemma:construction-structure-sheaf},
we obtain a spectrally ringed 
space $(|\Proj(\pi_0(A))|,\mathcal{O})$.
We shall show that 
it is a spectral scheme.

\begin{thm}
The spectrally ringed space
$(|\Proj(\pi_0(A))|,\mathcal{O})$
is a spectral scheme.
\end{thm}

\begin{proof}
In the proof of Lemma~\ref{lemma:construction-structure-sheaf},
we have showed that
the restriction
$(D_+(f),\mathcal{O}|_{D_+(f)})$
is equivalent to the affine spectral 
scheme $\Spec(R)$.
Since the set of all open subsets of the form $D_+(f)$
is a basis of the topology of $|\Proj(\pi_0(A))|$,
the theorem follows from
\cite{SAG}*{Corollary~1.1.6.4}.
\end{proof}

\begin{defn}
We denote  by 
\[ \Proj(A) \] 
the spectral scheme $(|\Proj(\pi_0(A))|,\mathcal{O})$
for a connective $\N$-graded $\E$-ring $A$.
We call $\Proj(A)$
the spectral projective scheme associated to $A$.
\end{defn}

\subsection{Properties of spectral projective schemes
associated to connective $\N$-graded $\E$-rings}
\label{subsection:properties-projective-schemes}

In this subsection we study
properties of spectral projective schemes
associated to connective $\N$-graded $\E$-rings. 
We show that Lurie's projective space
is obtained from a projective scheme for
some connective $\N$-graded $\E$-ring.

Let $A$ be a connective $\N$-graded $\E$-ring
and let $\Proj(A)=(|\Proj(\pi_0(A))|,\mathcal{O})$
be the spectral projective scheme associated to $A$.
We will show that there is a morphism
of spectral schemes from $\Proj(A)$
to $\Spec(A_0)$. 

\begin{construction}
The functor
$\mathcal{O}':\mathcal{U}(A)_e^{\mathrm{op}}\to\CAlg_{A_0}$
has been constructed 
in Construction~\ref{construction:spectral-projective-scheme},
which assigns to $D_+(f)$
the $\E$-ring $A[f^{-1}]_0$.
Let $\mathcal{O}_{A_0}: \mathcal{U}(A)^{\mathrm{op}}\to\CAlg_{A_0}$
be a right Kan extension of $\mathcal{O}'$
along the inclusion $\mathcal{U}(A)_e^{\mathrm{op}}\hookrightarrow
\mathcal{U}(A)^{\mathrm{op}}$.
Since the forgetful functor
$U: \CAlg_{A_0}\to \CAlg$
preserves small limits,
the structure sheaf $\mathcal{O}$ 
is equivalent to $U\circ \mathcal{O}_{A_0}$.
In particular,
we can regard $\mathcal{O}(|\Proj(\pi_0(A))|)$
as an $\E$-ring over $A_0$.
By \cite{SAG}*{Proposition~1.1.5.5},
we obtain a morphism
\[ \pi: \Proj(A)\longrightarrow \Spec(A_0) \]
of spectral schemes. 
\end{construction}

\begin{lemma}
The morphism $\pi: \Proj(A)\to \Spec(A_0)$
is separated.
\end{lemma}

\begin{proof}
Since the morphism $\Proj(\pi_0(A))\to \mathrm{\pi_0(A_0)}$
of the underlying schemes
is separated,
the lemma follows from
\cite{SAG}*{Remark~3.2.1.7}.
\end{proof}

\begin{prop}
Let $A$ be a connective $\N$-graded $\E$-ring.
Suppose that $\pi_0(A)$
is finitely generated as an $\N$-graded commutative ring
over $\pi_0(A_0)$.
Then the map $\pi: \Proj(A)\to\Spec(A_0)$
is proper.
\end{prop}

\begin{proof}
Since the morphism $\Proj(\pi_0(A))\to \mathrm{\pi_0(A_0)}$
of the underlying schemes
is proper,
the proposition follows from
\cite{SAG}*{Remark~5.1.2.2}.
\end{proof}

Let $\varphi: A\to B$ be a morphism of
connective $\N$-graded $\E$-rings.
Then $\varphi$ induces a map 
$\varphi_0: \pi_0(A)\to \pi_0(B)$
of $\N$-graded commutative rings.
We set $X=\Proj(A)$ and $Y=\Proj(B)$.
We denote by $\mathsf{X}=\Proj(\pi_0(A))$
and $\mathsf{Y}=\Proj(\pi_0(B))$
the underlying ordinary schemes of $X$ and $Y$, respectively.
We define $G(\varphi)$ to be an open subset of
$|\mathsf{Y}|$ given by
\[ 
G(\varphi)=\{\mathcal{P}\in |\mathsf{Y}||\ 
\mbox{$\varphi_0(\pi_0(A_i))\not\subset\mathcal{P}_i$ 
for some $i>0$} \}.
\]
Then $(G(\varphi),\mathcal{O}_Y|_{G(\varphi)})$
is an open subscheme of $Y$.
We have a map
\[ 
|\varphi^a|:G(\varphi)\to |\mathsf{X}| 
\]
of topological spaces
by assigning $\varphi_0^{-1}(\mathcal{P})$
to $\mathcal{P}\in G(\varphi)$.

\begin{prop}
There is an affine morphism
\[ 
\varphi^a: (G(\varphi),\mathcal{O}_Y|_{G(\varphi)})\to 
\Proj(A) 
\]
of spectral schemes
whose underlying map of topological spaces
is $|\varphi^a|$.
\end{prop}

\begin{proof}
Set $\mathcal{O}'=\mathcal{O}_Y|_{G(\varphi)}$.
We shall construct
a morphism $\mathcal{O}_X\to 
|\varphi^a|_*\mathcal{O}'$
of $\CAlg$-valued sheaves on $|\mathsf{X}|$.
Using the fact that 
$\CAlgzar_C(\Z)$
is equivalent
to an ordinary category
for a $\Z$-graded $\E$-ring $C$
by Theorem~\ref{thm:open-Zariski-immersions},
we can construct a natural transformation
$\mathcal{O}_X|_{\mathcal{U}(A)_e^{\mathrm{op}}}\to
(|\varphi^a|_*\mathcal{O}')|_{\mathcal{U}(A)_e^{\mathrm{op}}}$
given by
$A[f^{-1}]_0\to B[\varphi_0(f)^{-1}]_0$
on $D_+(f)$
for any homogeneous 
$f\in \pi_0(A)$ of positive degree.
Since $|\varphi^a|_*\mathcal{O}'$ is a right Kan
extension of its restriction to $\mathcal{U}(A)_e^{\rm op}$,
we obtain a morphism $\mathcal{O}_X\to |\varphi^a|_*\mathcal{O}'$
of $\CAlg$-valued sheaves on $|\mathsf{X}|$.
Since $|\varphi^a|^{-1}(D_+(f))=D_+(\varphi_0(f))$,
we see that $|\varphi^a|$ together with
$\mathcal{O}_X\to |\varphi^a|_*\mathcal{O}'$
defines an affine morphism
of spectral schemes. 
\end{proof}

Now, we suppose that
$\varphi: A\to B$ induces 
a surjection $\pi_0(A)\to \pi_0(B)$.
In this case we have $G(\varphi)=|\mathsf{Y}|$,
and hence we have an affine morphism
$\varphi^a: \Proj(B)\to\Proj(A)$. 

\begin{prop} 
Let $\varphi: A\to B$ be a morphism
of connective $\N$-graded $\E$-rings such that
$\varphi_0: \pi_0(A)\to \pi_0(B)$ is surjective.
Then the morphism
\[ 
\varphi^a: \Proj(B)\to \Proj(A) 
\]
is a closed immersion.
\end{prop}

\begin{proof}
For any homogeneous
$f\in\pi_0(A)$ of positive degree,
the restriction of $\varphi^a$ to ${D_+(\varphi_0(f))}$
is induced by the map
$A[f^{-1}]_0\to B[\varphi_0(f)^{-1}]_0$ of 
connective $\E$-rings.
Since this map induces a surjection 
$\pi_0(A[f^{-1}]_0)\to \pi_0(B[\varphi_0(f)^{-1}]_0)$
of commutative rings,
we see that $\varphi^a$ is a closed immersion
by \cite{SAG}*{Theorem~3.1.2.1}.
\end{proof}

In the remainder of this subsection
we will compare our construction of
spectral projective schemes with Lurie's projective spaces
defined in \cite{SAG}*{\S5.4.1}.
First, we recall the construction of Lurie's projective
spaces.

\begin{construction}[{cf.~\cite{SAG}*{Construction~5.4.1.3}}]
Let $R$ be an $\E$-ring.
For a $\Z$-graded set $X$, 
we denote by $R[X]$
the $\Z$-graded $R$-module given by $R\otimes\Sigma^{\infty}_+X$.
This construction $X\mapsto R[X]$
determines a symmetric monoidal functor
$\mathrm{Set}(\Z)\to \Mod_R(\Z)$,
where $\mathrm{Set}(\Z)$
is the category of $\Z$-graded sets.
Consequently,
this induces a functor
$\mathrm{CMon}(\Z)\to \CAlg_R(\Z)$,
where $\mathrm{CMon}(\Z)$
is the category of $\Z$-graded commutative monoids.

Let $[n]=\{0<1<\cdots<n\}$ and
let $P^{\circ}([n])$ be the set of all
nonempty subset of $[n]$.
For $I\in P^{\circ}([n])$,
we set
\[ M_I=\{(k_0,\ldots,k_n)\in\Z^{n+1}|\
\mbox{$k_0+\cdots k_n=0$ and $k_i\ge 0$ for $i\not\in I$}\}. \]
Then $M_I$ is a ($\{0\}$-graded) 
commutative monoid which depends
functorially on $I$.
Suppose that $R$ is a connective $\E$-ring.
By composing the functor 
$\Spet(-)$ with the above construction,
we obtain an affine spectral Deligne-Mumford
stack 
\[ \Spet(R[M_I]) \]
for each $I\in P^{\circ}([n])$
(see 
~\cite{SAG}*{\S1.4} 
for the definitions of $\Spet(-)$
and (affine) spectral Deligne-Mumford stacks). 
The construction
$I\mapsto \mbox{Sp\'{e}t}(R[M_I])$
determines a functor
\[ P^{\circ}([n])^{\mathrm{op}}\to \SpDM, \]
where $\SpDM$ is the $\infty$-category of
spectral Deligne-Mumford stacks.
Lurie's projective space 
\[ \mathbf{P}^n_R \]
is defined to be the colimit
$\colim_{I\in P^{\circ}([n])}\mbox{Sp\'{e}t}(R[M_I])$
in $\SpDM$.
\end{construction}

We set $N=\N^{n+1}$.
We regard $N$ as an $\N$-graded commutative monoid,
in which the degree $q$ component 
$N_q$ is the subset
\[ 
\{(k_0,\ldots,k_n)\in\N^{n+1}|\ k_0+\cdots+k_n=q\}. 
\]
Then we have a connective $\N$-graded $\E$-ring $R[N]$
and hence we obtain a spectral projective scheme 
\[ 
\Proj(R[N]). 
\]

By \cite{SAG}*{Remark~1.6.6.5},
there is a fully faithful functor 
\[ 
\iota: \mathrm{SpSch}\to\mathrm{SpDM}, 
\]
where $\mathrm{SpSch}$ is the $\infty$-category
of spectral schemes.
A spectral Deligne-Mumford stack is said to be
schematic if it belongs to the essential image of $\iota$. 
We shall show that the schematic spectral Deligne-Mumford
stack $\iota(\Proj(R[N]))$ is equivalent
to Lurie's projective space $\mathbf{P}^n_R$.

\begin{prop}
The schematic spectral Deligne-Mumford stack
$\iota(\Proj(R[N]))$
associated to the connective $\N$-graded $\E$-ring
$R[N]$ 
is equivalent to Lurie's projective space
$\mathbf{P}^n_R$.
\end{prop}

\begin{proof}
Note that 
if $I\subset J$ for $I,J\in P^{\circ}([n]))$,
then the map $\Spec(R[M_I])\to\Spec(R[M_J])$
is an open immersion.
We have 
$\Proj(R[N])\simeq
\colim_{I\in P^{\circ}([n])}\Spec(R[M_I])$.
Since $\iota(\Spec(A))\simeq 
\Spet(A)$ for any $\E$-ring $A$,
we see that 
$\iota(\Proj(R[N]))\simeq
\colim_{I\in P^{\circ}([n])}\mbox{Sp\'{e}t}(R[M_I])
\simeq \mathbf{P}^n_R$.
\end{proof}

\section{Quasi-coherent sheaves on projective spectral schemes}
\label{section:quasi-coheren-sheaf}

In this section we study the $\infty$-category
of quasi-coherent sheaves on 
spectral projective schemes.
First, we introduce a notion of Noetherian rings 
in the setting of $\N$-graded $\E$-rings.

\begin{defn}
Let $A$ be a connective $\N$-graded $\E$-ring.
We say that $A$ is Noetherian 
if $\pi_0(A)$ is a Noetherian $\N$-graded commutative ring
and $\pi_n(A)$ is a finitely generated $\N$-graded 
$\pi_0(A)$-module for any $n\in\Z$.
\end{defn}


Let $A$ be a connective Noetherian 
$\N$-graded $\E$-ring. 
In the following of this section
we assume that
$A$ satisfies the following condition: 

\begin{condition}\label{condition:N-graded-E-infinity-ring}\rm
There are finitely many elements
of $\pi_0(A_1)$
which generate $\pi_0(A)$
as an $\N$-graded commutative ring over $\pi_0(A_0)$.
\end{condition}

Let $X$ be the spectral projective scheme
$\Proj(A)$ associated to a connective Noetherian
$\E$-ring $A$ satisfying 
Condition~\ref{condition:N-graded-E-infinity-ring}.
In \S\ref{subsection:functor_widetilde}
we construct 
a symmetric monoidal functor
$\widetilde{(-)}$ from the $\infty$-category of 
$\Z$-graded $A$-modules to
the $\infty$-category of quasi-coherent sheaves on
$X$.
In \S\ref{subsection:right-adjoint-gamma}
we show that $\widetilde{(-)}$ admits a
right adjoint functor $\Gamma_*(X,-)$,
which is a graded extension of the global sections functor.
Finally,
we show that the functor $\Gamma_*(X,-)$ is fully faithful
and hence $\widetilde{(-)}$ is a localization functor.

\subsection{Construction of the functor $\widetilde{(-)}$}
\label{subsection:functor_widetilde}

In this subsection we shall construct a functor $\widetilde{(-)}$
from the $\infty$-category
of $\Z$-graded $A$-modules to
the $\infty$-category of quasi-coherent sheaves
on $\Proj(A)$.

Let $A$ be a connective Noetherian $\N$-graded $\E$-ring 
satisfying Condition~\ref{condition:N-graded-E-infinity-ring}.
Set $X=\Proj(A)$.
We have the $\infty$-category
$\QCoh(X)$ of quasi-coherent sheaves of $\mathcal{O}_X$-modules
on $X$ (cf.~\cite{SAG}*{Definition~2.2.2.1}).
First, we shall recall a description
of $\QCoh(X)$ in terms of $\infty$-categories
of modules over $\E$-rings.

We take a set $\{a_i\}_{i=1}^r\subset \pi_0(A_1)$ of generators 
of $\pi_0(A)$ as an $\N$-graded commutative
ring over $\pi_0(A_0)$.
We define a $\Z$-graded $\E$-ring $B$ by 
\[ B=A[a_1^{-1}]\times\cdots\times A[a_r^{-1}].\]
Let 
\[ g: A\to B \]
be the canonical morphism
of $\Z$-graded $\E$-rings.
We take a \v{C}ech nerve 
\[ 
C(g)^{\bullet}_+
\]
of $g$ 
in the opposite $\infty$-category of $\CAlg(\Z)$.
Then $C(g)^{\bullet}_+$ is an augmented cosimplicial object
of $\CAlg(\Z)$ such that 
$C(g)^{-1}_+\simeq A$ and $C(g)^n_+\simeq B^n$ for $n\ge 0$,
where $B^n$ is given by
\[ B^n=\overbrace{B\otimes_A \cdots\otimes_A B}^{n+1} .\]
Recall that we
have the functor $(-)_0: \CAlg(\Z)\to \CAlg$
in Definition~\ref{definition:restriction-CAlg}.
We denote by 
\[ C(g)^{\bullet}_0 \] 
the composite of $(-)_0$ 
with the restriction
$C(g)^{\bullet}=C(g)^{\bullet}_+|_{\Delta}$.

Let $\CAlg^{\mathrm{cn}}$ 
be the $\infty$-category of connective $\E$-rings
and 
let $\widehat{\mathcal{S}}$ be 
the very large $\infty$-category
of large spaces.
To a spectral scheme $X$,
we can associate a functor
\[ h_X: \CAlg^{\mathrm{cn}}\longrightarrow \widehat{\mathcal{S}} \] 
given by
\[ h_X(R)\simeq \Map_{\mathrm{SpSch}}(\Spec(R),X)\]
for $R\in\CAlg^{\mathrm{cn}}$.
By \cite{SAG}*{Proposition~1.6.3.3(2)},
we have a fully faithful functor
\[ h: \mathrm{SpSch}\longrightarrow 
      \Fun(\CAlg^{\mathrm{cn}},\widehat{\mathcal{S})}\]
which assign to $X\in\mathrm{SpSch}$ 
the functor $h_X$.
Let 
\[ \widehat{\mathcal{S}\mathrm{hv}}_{\mathrm{fpqc}} \]
be the full subcategory of 
$\Fun(\CAlg^{\mathrm{cn}},\widehat{\mathcal{S}})$
spanned by those functors which are sheaves 
for the fpqc topology.
By \cite{SAG}*{Theorem~1.6.2.1},
we can regard a spectral scheme $X$ as an object of 
$\widehat{\mathcal{S}\mathrm{hv}}_{\mathrm{fpqc}}$
through the functor $h$.

There is a faithfully flat affine morphism
\[ f: U\longrightarrow X, \] 
where $U={\rm Spec}(B_0)$.
We denote by $U_{\bullet}^+$ the \v{C}ech nerve of $f$
in the $\infty$-category of spectral schemes.
Then we have 
$U_{-1}^+\simeq X$
and
$U_n^+\simeq {\rm Spec}(B^n_0)$
for $n\ge 0$.
We denote by
\[ U_{\bullet}\]
the simplicial object 
obtained by the restriction of $U^+_{\bullet}$ to ${\Delta}$.
Note that there is an equivalence
\[ U_{\bullet}\simeq {\rm Spec}(C(g)^{\bullet}_0) \]
of simplicial objects of affine spectral schemes.
Since $f: U\to X$ is an effective epimorphism,
we have an equivalence
\[ |U_{\bullet}|\stackrel{\simeq}{\longrightarrow}
    X \]
in 
$\widehat{\mathcal{S}\mathrm{hv}}_{\mathrm{fpqc}}$,
where the left hand side is the geometric realization
of the simplicial object $U_{\bullet}$.

By \cite{SAG}*{\S2.2.2},
we have an $\infty$-category $\QCoh(X)$
of quasi-coherent sheaves of $\mathcal{O}_X$-modules on $X$.
The $\infty$-category $\QCoh(X)$
is stable and presentable by \cite{SAG}*{Proposition~2.2.4.1}.
Furthermore,
$\QCoh(X)$ admits a symmetric monoidal structure
with unit $\mathcal{O}_X$
by \cite{SAG}*{Proposition~2.2.4.2}.
The following lemma describes $\QCoh(X)$
in terms of $\infty$-categories of modules
over $\E$-rings.

\begin{lemma}\label{lemma:equivalence-QCohX-lim-Mod-B}
There is an equivalence
\[ \QCoh(X)\simeq \lim_{\Delta} \Mod_{B^{\bullet}_0} \]
of symmetric monoidal stable
$\infty$-categories.
\end{lemma}

\begin{proof}
By \cite{SAG}*{Proposition~6.2.3.1(b)},
there is an equivalence
\[ \QCoh(X)\stackrel{\simeq}{\longrightarrow}
   \lim_{\Delta} \QCoh(U_{\bullet})\]
of stable $\infty$-categories.
Since $\QCoh(\Spec R) \simeq \Mod_R$ 
for any $\E$-ring $R$
by \cite{SAG}*{Proposition 2.2.3.3.},
we obtain an equivalence
\[ \QCoh(X)\stackrel{\simeq}{\longrightarrow}
   \lim_{\Delta} \Mod_{B^{\bullet}_0} \]
of stable $\infty$-categories.
By \cite{SAG}*{\S6.2.6},
we see that 
$\QCoh(X)\simeq \lim_{\Delta} \Mod_{B^{\bullet}_0}$
is an equivalence of symmetric monoidal
$\infty$-categories.
\end{proof}

Next, we shall construct a functor
from $\Mod_A(\Z)$ to $\QCoh(X)$.
Recall that $\Sp(\Z)$
is a symmetric monoidal stable presentable $\infty$-category,
in which the tensor product commutes with small colimits
separately in each variable
by Propositions~\ref{prop:graded-spectra-stable-presentable} and
\ref{prop:symmetric-monoidal-structure-graded-spectra}.
In other words,
$\Sp(\Z)$ is a commutative algebra object
of ${\mathcal Pr}^{\mathrm{L}}_{\mathrm{St}}$,
where  ${\mathcal Pr}^{\mathrm{L}}_{\mathrm{St}}$ is 
the $\infty$-category of stable presentable
$\infty$-categories and 
colimit-preserving functors.
We denote by
$\CAlg({\mathcal Pr}^{\mathrm{L}}_{\mathrm{St}})$
the $\infty$-category of commutative algebra objects
of ${\mathcal Pr}^{\mathrm{L}}_{\mathrm{St}}$.
By \cite{HA}*{Theorem~4.8.5.16},
we have a functor
\[ \Mod_{(-)}(\Z): \CAlg(\Z)\longrightarrow 
   \CAlg({\mathcal Pr}^{\mathrm{L}}_{\mathrm{St}}),\]
which assigns to $D\in\CAlg(\Z)$
the symmetric monoidal $\infty$-category
$\Mod_D(\Z)$ of $\Z$-graded $D$-modules.
By applying the functor $\Mod_{(-)}(\Z)$
to $C(g)^{\bullet}_+$
and using the equivalences
$\Mod_{B^n}(\Z)\simeq\Mod_{B^n_0}$
for $n\ge 0$
by Corollary~\ref{cor:periodic-equivalence},
we obtain 
a symmetric monoidal functor
\[ \Mod_A(\Z)\longrightarrow
   \lim_{\Delta}\Mod_{B_0^{\bullet}}.\]

\begin{defn}
We define a functor
\[ \widetilde{(-)}: \Mod_A(\Z)\longrightarrow
   \QCoh(X) \]
to be the composite of the functor
$\Mod_A(\Z)\to\lim_{\Delta}\Mod_{B_0^{\bullet}}$
with the equivalence between
$\lim_{\Delta}\Mod_{B_0^{\bullet}}$ and 
$\QCoh(X)$
by Lemma~\ref{lemma:equivalence-QCohX-lim-Mod-B}.
We call $\widetilde{M}$
the quasi-coherent sheaf on $X$ associated to
a $\Z$-graded $A$-modules $M$.
\end{defn}

By the construction,
the functor $\widetilde{(-)}: \Mod_{A}(\Z)\to
\QCoh(X)$ is 
symmetric monoidal.
In particular, we have an equivalences
\[ \mathcal{O}_X\simeq \widetilde{A}\]
and 
\[ \widetilde{M}\otimes_{\mathcal{O}_X}\widetilde{N}
   \simeq \widetilde{M\otimes_AN}\]
for any $M,N\in\Mod_A(\Z)$.

Recall that we have defined the shifting functor
$(q): \Mod_A(\Z)\to\Mod_A(\Z)$
in Definition~\ref{definition:twisting-module-spectra},
which is given by $M(q)_n\simeq M_{q+n}$
for $M\in\Mod_A(Z)$ and $q,n\in\Z$.

\begin{defn}
For $q\in\Z$,
we define 
\[ \mathcal{O}_X(q) \] 
to be the quasi-coherent sheaf $\widetilde{A(q)}$ on $X$.
\end{defn}

Since 
$\mathcal{O}_X\simeq \widetilde{A}$
and
$\widetilde{M}\otimes_{\mathcal{O}_X}\widetilde{N}\simeq
\widetilde{M\otimes_A N}$,
we have 
$\mathcal{O}_X(0)\simeq\mathcal{O}_X$
and
$\mathcal{O}_X(q) \otimes_{\mathcal{O}_X}
\mathcal{O}_X(q') \simeq \mathcal{O}_X(q + q')$
for any $q,q'\in\Z$.
By \cite{SAG}*{Proposition~2.9.4.2},
we see that $\mathcal{O}_X(q)$ is a line bundle
on $X$ for any $q\in\Z$, that is,
it is locally free of rank $1$.

\begin{defn}
For a quasi-coherent sheaf $\mathscr{F}$
of $\mathcal{O}_X$-modules on $X$ and $q\in\Z$,
we define 
\[ \mathscr{F}(q) \]
to be the tensor product
$\mathscr{F} \otimes_{\mathcal{O}_X} 
\mathcal{O}_X (q)$.
\end{defn}

Note that if $\mathscr{F}$ is the quasi-coherent
sheaf associated to a $\Z$-graded $A$-module 
$M$,
then we have
an equivalence 
$\mathscr{F}(q)\simeq \widetilde{M(q)}$
of quasi-coherent sheaves.

\subsection{The functor $\Gamma_*(X,-)$}
\label{subsection:right-adjoint-gamma}

We have constructed a functor
$\widetilde{(-)}: \Mod_A(\Z)\to \QCoh(X)$
in \S\ref{subsection:functor_widetilde}.
In this subsection we shall construct
a functor $\Gamma_*(X,-): \QCoh(X)\to\Mod_A(\Z)$,
which is a right adjoint to the functor $\widetilde{(-)}$.
We show that $\Gamma_*(X,-)$ is fully faithful
and hence $\widetilde{(-)}$ is a localization functor.

In \S\ref{subsection:properties-projective-schemes}
we have constructed a morphism of spectral schemes
$\pi: X\longrightarrow {\rm Spec}(A_0)$,
where $A_0$ is the $\E$-ring given by
the degree $0$ part of $A$.
Since we have an equivalence $\QCoh({\rm Spec}(A_0))\simeq
\Mod_{A_0}$,
we obtain a pullback functor
\[ \pi^*: \Mod_{A_0}\longrightarrow \QCoh(X), \]
which assigns to $N\in\Mod_{A_0}$
the quasi-coherent sheaf on $X$
associated to the $\Z$-graded $A$-module
$A\otimes_{A_0}N$.
The pushforward functor 
\[ \pi_*: \QCoh(X)\longrightarrow \Mod_{A_0} \]
is a right adjoint to $\pi^*$.
Let $\Upsilon: \Mod_{A_0}\to{\rm Sp}$
be the forgetful functor.
Then $\Upsilon\circ\pi_*$ is equivalent
to the global sections functor $\Gamma(X,-)$.

We consider a right adjoint 
to the functor 
$\widetilde{(-)}: \Mod_A(\Z)\to \QCoh(X)$.

\begin{prop}
There is a right adjoint
\[ \Gamma_*(X,-): \QCoh(X)\longrightarrow \Mod_A(\Z) \]
to the functor $\widetilde(-)$.
\end{prop}

\begin{proof}
Recall that the functor $\widetilde{(-)}$
is obtained from the augmented
cosimplicial diagram
$\Mod_{(-)}(\Z)\circ C(g)^{\bullet}_+:
\Delta_+\to 
\CAlg(\Z)\to \widehat{\Cat{\infty}}$.
Since the functor
$\Mod_{(-)}(\Z): \CAlg(\Z)\to \widehat{\Cat{\infty}}$
factors through the $\infty$-category 
${\mathcal Pr}^{\mathrm{L}}_{\mathrm{St}}$
of stable presentable $\infty$-categories and
colimit-preserving functors,
we see that the functor
$\widetilde{(-)}: \Mod_A(\Z)\to\QCoh(X)$
is a morphism in ${\mathcal Pr}^{\mathrm{L}}_{\mathrm{St}}$, 
and hence there exists a right adjoint
$\Gamma_*(X,-): \QCoh(X)\longrightarrow \Mod_A(\Z)$
to $\widetilde{(-)}$.
\end{proof}

Now we shall describe the functor $\Gamma_*(X,-)$
more explicitly.
Let 
\[ p: E\to \Delta_+ \]
be the coCartesian fibration
associated to the functor
$\Mod_{(-)}(\Z)\circ C(g)^{\bullet}_+:
\Delta_+\to\widehat{\Cat{\infty}}$.
We have equivalences $E_{-1}\simeq\Mod_A(\Z)$
and $E_{n}\simeq\Mod_{B^n}(\Z)$ for $n\ge 0$. 

Let $\mathscr{F}$ be a quasi-coherent sheaf of
$\mathcal{O}_X$-modules on $X$.
For $n\ge 0$,
we let 
\[ \Gamma(U_n,\mathscr{F}(*)) \]
be the $\Z$-graded $B^n$-module,
in which
the degree $q$-component is given by
$\Gamma(U_n,f_n^*\mathscr{F}(q))$,
where $f_n: U_n\to X$ is the canonical map.
By \cite{HT}*{Proposition~3.3.3.2},
the quasi-coherent sheaf $\mathscr{F}$ 
can be identified with
a coCartesian section $s$
of $p$ over $\Delta\subset\Delta_+$:
\[ \xymatrix{
     & E \ar[d]^p\\
    \Delta\ar@{..>}^s[ur]\ar@{^(_->}[r]&\Delta_+
   }\]
in which
$s([n])\simeq \Gamma(U_n,\mathscr{F}(*))$ for $[n]\in\Delta$.
Since $\Mod_{(-)}(\Z)\circ C(g)^{\bullet}_+$
factors through ${\mathcal Pr}^{\mathrm{L}}_{\mathrm{St}}$,
the map $p: E\to\Delta_+$ is a presentable fibration
by \cite{HT}*{Proposition~5.5.3.3(2)}.
In particular, $p: E\to\Delta_+$ is a Cartesian fibration.
For the coCartesian section $s$,
by taking the Cartesian pullback,
we obtain a cosimplicial object
\[ \Gamma(U_{\bullet},\mathscr{F}(*)):
   \Delta\longrightarrow \Mod_A(\Z), \]
which is given by
$[n]\mapsto \Gamma(U_n,\mathscr{F}(*))$,
where we regard $\Gamma(U_n,\mathscr{F}(*))$
as a $\Z$-graded $A$-module 
through the canonical map $A\to B^n$.

\begin{prop}
For any quasi-coherent sheaf $\mathscr{F}$ of
$\mathcal{O}_X$-modules on $X$,
we have an equivalence
\[ \Gamma_*(X,\mathscr{F})\simeq
   \lim_{\Delta}\Gamma(U_{\bullet},\mathscr{F}(*))\]
of $\Z$-graded $A$-modules.
\end{prop}

\begin{proof}
Recall that we have the equivalence
$\QCoh(X)\stackrel{\simeq}{\longrightarrow}
   \lim_{\Delta} \Mod_{B^{\bullet}}(\Z)$
of $\infty$-categories.
This implies an equivalence
\[ \Map_{\QCoh(X)}(\widetilde{M},\mathscr{F})\simeq
   \lim_{\Delta}\Map_{\Mod_B^{\bullet}(\Z)}
   (B^{\bullet}\otimes_AM,\Gamma(U_{\bullet},\mathscr{F}(*))) \]
of mapping spaces.
Since we have an equivalence
\[ \Map_{\Mod_B^{\bullet}(\Z)}
   (B^{\bullet}\otimes_AM,\Gamma(U_{\bullet},\mathscr{F}(*)))
   \simeq   
   \Map_{\Mod_A(\Z)}
   (M,\Gamma(U_{\bullet},\mathscr{F}(*)))\] 
of cosimplicial spaces,
there is a natural equivalence
\[ \Map_{\QCoh(X)}(\widetilde{M},\mathscr{F})\simeq
   \Map_{\Mod_A(\Z)}
   (M,\lim_{\Delta}\Gamma(U_{\bullet},\mathscr{F}(*))). \]
Hence we obtain
$\Gamma_*(X,\mathscr{F})\simeq
\lim_{\Delta}\Gamma(U_{\bullet},\mathscr{F}(*))$.
This completes the proof.
\end{proof}

Next,
we shall show that the functor
$\Gamma_*(X,-)$ is fully faithful.
Recall that $B=A[a_1^{-1}]\times\cdots A[a_r^{-1}]$,
where $\{a_i\}_{i=1}^r\subset \pi_0(A_1)$ is 
the set of generators 
of $\pi_0(A)$ as an $\N$-graded commutative
ring over $\pi_0(A_0)$.
We have the faithfully flat affine morphism
$f: U\to X$, where $U={\rm Spec}(B_0)$.
Note that there is an equivalence
$\Gamma(U,f^*\mathcal{O}_X(*))\simeq B$
of $\mathbb{Z}$-graded $\mathbb{E}_{\infty}$-rings
and hence that $\Gamma(U,f^*\mathscr{F}(*))$
is a $\mathbb{Z}$-graded $B$-module
for a quasi-coherent sheaf $\mathscr{F}$
of $\mathcal{O}_X$-modules on $X$.
We have the restriction map 
$\Gamma_*(X,\mathscr{F})\to \Gamma(U,f^*\mathscr{F}(*))$,
which induces a map
\[ B\otimes_A\Gamma_*(X,\mathscr{F})\to \Gamma(U,f^*\mathscr{F}(*)) \]
of $\Z$-graded $B$-modules.
We shall compare $\Gamma(U,f^*\mathscr{F}(*))$
with the scalar extension 
$B\otimes_A\Gamma_*(X,\mathscr{F})$.

\begin{lemma}\label{lemma:compatibility-of-localization}
Let $\mathscr{F}$ be a quasi-coherent sheaf 
of $\mathcal{O}_X$-module on $X$.
There is a natural equivalence
\[ B\otimes_A\Gamma_*(X,\mathscr{F}) 
   \stackrel{\simeq}{\longrightarrow}  
   \Gamma(U,f^*\mathscr{F}(*))  \]
of $\mathbb{Z}$-graded $B$-modules.
\end{lemma}

\begin{proof}
We have $U=V_1\times\cdots \times V_r$,
where $V_i={\rm Spec}(A[a_i^{-1}]_0)$
for $1\le i\le r$.
This implies a decomposition
\[ \Gamma(U,f^*\mathscr{F}(*))\simeq
   \Gamma(V_1,\mathscr{F}(*))\times\cdots\times
   \Gamma(V_r,\mathscr{F}(*)),\]
where $\Gamma(V_i,\mathscr{F}(*))$
is a $\Z$-graded $A[a_i^{-1}]$-module
for $1\le i\le r$. 
Since $\Gamma(V_i,\mathscr{F}(*))$
is a $\mathbb{Z}$-graded $A[a_i^{-1}]$-module,
the restriction map
$\Gamma_*(X,\mathscr{F})\to\Gamma(U,f^*\mathscr{F}(*))$
induces a map
$\Gamma_*(X,\mathscr{F})[a_i^{-1}]\to \Gamma_*(V_i,\mathscr{F})$
of $\mathbb{Z}$-graded $A[a_i^{-1}]$-modules.
It suffices to show that this map
is an equivalence for any $i$ with $1\le i\le r$.

Let $P$ be the partially ordered set of all nonempty
finite subsets of $\{1,\ldots,r\}$.
We set $V_I=\cap_{i\in I}V_i$ for $I\in P$. 
By \cite{SAG}*{Proposition 1.1.4.4}, 
we have an equivalence
\[ \Gamma_*(X,\mathscr{F})\simeq\
   \subrel{I\in P}{\lim}
   \Gamma(V_I,\mathscr{F}(*))\]
of $\mathbb{Z}$-graded $A$-modules.
Note that the right hand side is a finite
limit indexed by $P$.
Since filtered colimits commute 
with finite limits,
we obtain an equivalence
\[ \Gamma_*(X,\mathscr{F})[a_i^{-1}]\simeq\
   \subrel{I\in P}{\lim}
   (\Gamma(V_I,\mathscr{F}(*))[a_i^{-1}])\]
of $\mathbb{Z}$-graded $A[a_i^{-1}]$-modules.
By definition, we have an equivalence
\[ \Gamma(V_I,\mathscr{F}(*))[a_i^{-1}]\simeq
   \Gamma(V_{I\cup\{i\}},\mathscr{F}(*))\]
for any $I\in P$.
We consider a functor
$\theta: P\to {\rm Mod}_{A[a_i^{-1}]}(\Z)$ which assigns
to $I\in P$ the $\mathbb{Z}$-graded $A[a_i^{-1}]$-module 
$\Gamma(V_{I\cup \{i\}},\mathscr{F}(*))$.
Let $P_i$ be the subset of $P$
consisting of finite subsets of $\{1,\ldots,r\}$
which contain $i$.
Since the functor $\theta$ is a right Kan extension
of the restriction to $P_i$,
we have an equivalence
\[ \subrel{I\in P}{\lim}
   \Gamma(V_{I\cup\{i\}},\mathscr{F}(*))\simeq\
   \subrel{J\in P_i}{\lim}
   \Gamma(V_J,\mathscr{F}(*)).\]
By \cite{SAG}*{Proposition 1.1.4.4}, 
we have an equivalence
\[ \Gamma(V_i,\mathscr{F}(*))\simeq\
   \subrel{J\in P_i}{\lim}
   \Gamma(V_J,\mathscr{F}(*)).\]
Hence we obtain an equivalence
$\Gamma(V_i,\mathscr{F}(*))\simeq
   \Gamma_*(X,\mathscr{F})[a_i^{-1}]$
of $\mathbb{Z}$-graded $A[a_i^{-1}]$-modules.   
\end{proof}

\begin{prop}
The functor $\Gamma_*(X,-): \QCoh(X)\to \Mod_A(\Z)$
is fully faithful.
\end{prop}

\proof
We shall show that
the counit $\widetilde{(-)}\circ \Gamma_*(X,-)\to 
{\rm id}_{\QCoh(X)}$ of the adjunction is an equivalence.
Since the pullback functor 
$f^*: \QCoh(X)\to \QCoh(U)\simeq \Mod_{B_0}$
is conservative and
there is an equivalence
$B\otimes_{B_0}(-):\Mod_{B_0}\simeq \Mod_B(\Z)$
by Corollary~\ref{cor:periodic-equivalence},
it suffices to show that 
the map   
$B\otimes_{B_0}f^*\widetilde{\Gamma_*(X,\mathscr{F})}\to
B\otimes_{B_0}f^*\mathscr{F}$
is an equivalence
for any quasi-coherent sheaf $\mathscr{F}$ on $X$.
Note that the composite
$\Mod_A(\Z)\stackrel{\widetilde{(-)}}{\to} 
\QCoh(X)\stackrel{f^*}{\to}
\Mod_{B_0}\simeq \Mod_B(\Z)$
is given by
$M\mapsto B\otimes_AM$
and that the composite
$\QCoh(X)\stackrel{f^*}{\to}\Mod_{B_0}
\stackrel{\simeq}{\to}\Mod_{B}(\Z)$
is given by
$\mathscr{F}\mapsto
\Gamma(U,f^*\mathscr{F}(*))$.
Thus, we can identify the map
$B\otimes_{B_0}f^*\widetilde{\Gamma_*(X,\mathscr{F})}\to
B\otimes_{B_0}f^*\mathscr{F}$
with the map
$B\otimes_A\Gamma_*(X,\mathscr{F})\to
   \Gamma(U,f^*\mathscr{F}(*))$,
which is an equivalence 
by Lemma~\ref{lemma:compatibility-of-localization}.
This completes the proof.
\qed

\bigskip

We have the adjunction
\[ \widetilde{(-)}: \Mod_A(\Z)\rightleftarrows
                    \QCoh(X): \Gamma_*(X,-)\]
of $\infty$-categories.
Since the left adjoint $\widetilde{(-)}$ is 
a symmetric monoidal functor,
the right adjoint $\Gamma_*(X,-)$
is a lax symmetric monoidal functor.
In particular,
$\Gamma_*(X,\mathcal{O}_X)$
is a $\Z$-graded $\E$-ring and  
there is a map
\[ A\longrightarrow \Gamma_*(X,\mathcal{O}_X) \]
of $\Z$-graded $\E$-rings.
For a quasi-coherent sheaf $\mathscr{F}$
of $\mathcal{O}_X$-modules on $X$,
$\Gamma_*(X,\mathscr{F})$ is a
$\Z$-graded $\Gamma_*(X,\mathcal{O}_X)$-module.
We note that the $\Z$-graded $A$-module structure
on $\Gamma_*(X,\mathscr{F})$
is obtained from
the $\Z$-graded $\Gamma_*(X,\mathcal{O}_X)$-module structure
through the map
$A\to\Gamma_*(X,\mathcal{O}_X)$
of $\Z$-graded $\E$-rings.

\section{Almost perfect quasi-coherent sheaves on 
projective spectral schemes}
\label{section:almost_perfect_projective_scheme}

The Serre theorem describes the category of 
coherent sheaves 
on an ordinary projective scheme 
$\Proj(R)$ 
associated to an $\N$-graded commutative ring $R$ as
the quotient of the abelian category of
finitely generated $\Z$-graded $R$-modules by the Serre 
subcategory of torsion $\Z$-graded $R$-modules.
In this section we consider an analogue of the Serre theorem
in spectral algebraic geometry.
In \S\ref{subsection:ordinary_projective_scheme}
we recall some properties of the cohomology of coherent sheaves 
on ordinary projective schemes.
In \S\ref{subsection:almost_perfect_quasi-cohenre_sheaf}
we give a characterization
of the $\Z$-graded modules of global sections
for almost perfect quasi-coherent sheaves
on spectral projective schemes.
In \S\ref{subsection:Serre_theorem_spectral}
we prove the Serre theorem in spectral algebraic geometry.
We show that 
the $\infty$-category of almost perfect quasi-coherent
sheaves 
on a spectral projective scheme 
$\Proj(A)$ 
associated to a connective $\N$-graded $\E$-ring $A$ is 
equivalent to the $\infty$-category obtained 
as the Verdier quotient of 
the $\infty$-category of almost finitely generated 
$\Z$-graded $A$-modules 
by the full subcategory of almost torsion
$\Z$-graded $A$-modules.

\subsection{Ordinary projective schemes and
quasi-coherent sheaves over them}
\label{subsection:ordinary_projective_scheme}

In this subsection
we recall some properties of the cohomology of coherent sheaves 
on ordinary projective schemes
associated to $\N$-graded commutative rings.

First, we fix some terminology and notation on
ordinary graded modules.
Let $R$ be an $\N$-graded commutative ring.
We denote by 
$\Mod_R^{\heartsuit}(\Z)$
the category of $\mathbb{Z}$-graded $R$-modules
and $R$-module homomorphisms of degree $0$.
For a $\mathbb{Z}$-graded $R$-module $M$,
we denote by $M_n$ the degree $n$ component of $M$.
We say that $M$ is bounded above
if there exists $n_0 \in \mathbb{Z}$ such that
$M_n=0$ for $n>n_0$,
and that 
$M$ is locally bounded above
if any finitely generated $\mathbb{Z}$-graded $R$-submodule of $M$ 
is bounded above.

For $k\in\mathbb{Z}$, we have a functor
\[ \sigma_{\ge k}:
\Mod_R^{\heartsuit}(\Z)
\longrightarrow 
\Mod_R^{\heartsuit}(\Z), 
\]
which assigns to $M\in\Mod_R^{\heartsuit}(\Z)$ 
the $\mathbb{Z}$-graded $R$-module 
$\sigma_{\ge k}M$ given by
\[ (\sigma_{\ge k}M)_n=\left\{
               \begin{array}{cl}
                 M_n & (n\ge k),\\[2mm]
                 0   & (n<k).\\
               \end{array} \right.\]
There is a natural transformation
\[ \sigma_{\ge k}\longrightarrow 
   \mathrm{id}_{\Mod_R^{\heartsuit}(\Z)} \]
of functors such that $(\sigma_{\ge k}M)_n\to M_n$
is the identity map for each $n\ge k$.

\begin{defn}
We say that a $\Z$-graded $R$-module $M$ 
is strongly quasi-finitely generated if 
$\sigma_{\ge k}M$ is finitely generated 
as a $\Z$-graded $R$-module
for all $k\in \Z$. 
We define 
\[ \Mod_R^{\heartsuit\mathrm{sqfg}}(\Z) \]
to be the full subcategory of 
$\Mod_R^{\heartsuit}(\Z)$
consisting of strongly quasi-finitely generated 
$\Z$-graded $R$-modules.
\end{defn}

\begin{rem}
If a $\Z$-graded $R$-module $M$
is strongly quasi-finitely generated and locally bounded above,
then $M$ is bounded above.
\end{rem}

We recall that a nonempty full subcategory 
of an abelian category is said to be a Serre subcategory
if it is closed under extensions and subquotients. 
We easily obtain the following lemma.

\begin{lemma}\label{lemma:mod-sqfg-Serre-subcategory}
If $R$ is a Noetherian $\N$-graded commutative ring, 
then 
$\Mod_R^{\heartsuit\mathrm{sqfg}}(\Z)$
is a Serre subcategory of 
$\Mod_R^{\heartsuit}(\Z)$.
\end{lemma}

For an $\N$-graded commutative ring $R$,
we denote by $\Proj(R)$
the ordinary projective scheme associated to $R$.
In the following of this subsection
we assume that $R_0$ is Noetherian
and that $R$ is generated by
finitely many elements of $R_1$
as an $\N$-graded commutative ring over $R_0$.
Note that in this case $R$ is a Noetherian
$\N$-graded commutative ring.
We set $\mathsf{X}=\Proj(R)$.
For $M\in\Mod_R^{\heartsuit}(\Z)$,
we have a quasi-coherent sheaf 
$\widetilde{M}$ on $\mathsf{X}$ associated to $M$.
This construction induces a functor
\[ \widetilde{(-)}: 
   \Mod_R^{\heartsuit}(\Z)
   \longrightarrow
   \QCoh(\mathsf{X}), \]    
where $\QCoh(\mathsf{X})$ is the category
of quasi-coherent sheaves of $\mathcal{O}_{\mathsf{X}}$-modules
on $\mathsf{X}$.
Note that 
$\widetilde{M}= 0$ if and only if
$M$ is locally bounded above for 
$M\in \Mod_R^{\heartsuit}(\Z)$.
We denote by
$\Mod_R^{\heartsuit\fg}(\Z)$
the full subcategory of 
$\Mod_R^{\heartsuit}(\Z)$
consisting of finitely generated $\Z$-graded $R$-modules.
If $M\in \Mod_R^{\heartsuit\fg}(\Z)$,
then $\widetilde{M}$ is a coherent sheaf
of $\mathcal{O}_{\mathsf{X}}$-modules on $\mathsf{X}$. 
Thus, we obtain a functor
\[ \widetilde{(-)}: 
   \Mod_R^{\heartsuit\fg}(\Z)
   \longrightarrow
   \Coh(\mathsf{X})\]
by restricting the functor $\widetilde{(-)}$
to $\Mod_R^{\heartsuit\fg}(\Z)$,
where $\Coh(\mathsf{X})$ is 
the category of coherent sheaves of 
$\mathcal{O}_{\mathsf{X}}$-modules on $\mathsf{X}$.

For $M\in \Mod_R^{\heartsuit\fg}(\Z)$,
we have $\widetilde{M}= 0$ if and only if
$M$ is bounded above.
We say that a finitely generated $\Z$-graded $R$-module
is torsion if it is bounded above.
We denote by 
$\Mod_R^{\heartsuit\tor}(\Z)$
the full subcategory of 
$\Mod_R^{\heartsuit\fg}(\Z)$
consisting of torsion $\Z$-graded $R$-modules.
We note that 
$\Mod_R^{\heartsuit\tor}(\Z)$
is a Serre subcategory 
of the abelian category 
$\Mod_R^{\heartsuit\fg}(\Z)$.

The classical Serre theorem describes
the abelian category $\Coh(\mathsf{X})$
in terms of $\Z$-graded $R$-modules.

\begin{prop}[{cf.~\cite{Stack}*{Proposition~0BXD}}]
\label{prop:classical-Serre-theorem}
Let $R$ be an $\mathbb{N}$-graded commutative ring.
We assume that $R_0$ is Noetherian and that $R$ is generated by
finitely many elements of $R_1$ as an $\N$-graded 
commutative ring over $R_0$.
Set $\mathsf{X}=\Proj(R)$.
The functor $\widetilde{(-)}$ induces an equivalence
\[ \Mod_R^{\heartsuit\fg}(\Z)/\Mod_R^{\heartsuit\tor}(\Z)
   \stackrel{\simeq}{\longrightarrow}
   \Coh(\mathsf{X}) \]
of abelian categories,
where the left hand side is the quotient category of
the abelian category 
$\Mod_R^{\heartsuit\fg}(\Z)$
by the Serre subcategory 
$\Mod_R^{\heartsuit\tor}(\Z)$.
\end{prop}

Let $\mathcal{F}$ be a quasi-coherent sheaf of 
$\mathcal{O}_{\mathsf{X}}$-module on $\mathsf{X}$.
For each $p\ge 0$,
we have a $\Z$-graded $R$-module 
\[ H^p(\mathsf{X},\mathcal{F}(*)), \]
in which the degree $q$ component is given by
$H^p(\mathsf{X},\mathcal{F}(q))$.
Since $\mathsf{X}$ is a quasi-compact separated scheme,
there exists an integer $s_0$ such that
$H^s(\mathsf{X},\mathcal{F}(*))=0$ for 
any quasi-coherent sheaf $\mathcal{F}$
of $\mathcal{O}_X$-module on $\mathsf{X}$ and $s>s_0$.

On the cohomology of coherent sheaves
on ordinary projective schemes,
we have the following proposition.

\begin{prop}[{cf.~\cite{Stack}*{Lemma~0AG6}}]
\label{prop:cohomology_projective_scheme}
Let $R$ be an $\N$-graded commutative ring. 
We assume that $R_0$ is Noetherian and that $R$ is generated by
finitely many elements of $R_1$ as an $\N$-graded 
commutative ring over $R_0$.
If $\mathcal{F}$ is a coherent sheaf of $\mathcal{O}_{\mathsf{X}}$-modules
on $\mathsf{X}=\Proj(R)$,
then the $\mathbb{Z}$-graded $R$-module 
$H^p(\mathsf{X},\mathcal{F}(*))$
is strongly quasi-finitely generated
for all $p\ge 0$.
Furthermore,
$H^p(\mathsf{X},\mathcal{F}(*))$
is bounded above for $p>0$.
\end{prop}


\subsection{Almost perfect quasi-coherent sheaves
on spectral projective schemes}
\label{subsection:almost_perfect_quasi-cohenre_sheaf}

In this subsection we give a characterization
of $\Z$-graded modules of global sections
for almost perfect quasi-coherent sheaves
on spectral projective schemes 
associated to connective $\N$-graded $\E$-rings.

Let $A$ be an $\N$-graded $\mathbb{E}_{\infty}$-ring.
First, we introduce a notion of perfect modules 
in the setting of $\Z$-graded $A$-modules.

\begin{defn}
We let 
\[  \Mod_A^{\perf}(\Z) \]
be the smallest stable subcategory 
of $\Mod_A(\Z)$ which contains $A(q)$ for all $q \in\mathbb{Z}$
and is closed under retracts.
We say that a $\Z$-graded $A$-module $M$
is perfect
if it belongs to 
the full subcategory $\Mod_A^{\perf}(\Z)$.
\end{defn}

Now we introduce a notion of almost finitely generated
modules in the setting of $\Z$-graded $A$-modules.

\begin{defn}\label{afg}
Let $A$ be a connective $\N$-graded $\mathbb{E}_{\infty}$-ring.
We say that a $\Z$-graded $A$-module $M$ 
is almost finitely generated
if the following conditions are satisfied:
\begin{enumerate}[(i)]
\item For each $n\in\mathbb{Z}$,
the $\mathbb{Z}$-graded $\pi_0(A)$-module
$\pi_n(M)$ is finitely generated. 
\item For $n\ll 0$, $\pi_n(M)=0$.
\end{enumerate}
We define an $\infty$-category
\[ \Mod_A^{\afg}(\Z) \]
to be the full subcategory of $\Mod_A(\Z)$
spanned by almost finitely generated
$\mathbb{Z}$-graded $A$-modules. 
We denote by
\[ \Mod_A^{\afg}(\Z)_{\ge 0} \]
the full subcategory of $\Mod_A^{\afg}(\Z)$
consisting of almost finitely generated
$\Z$-graded $A$-modules $M$
such that $\pi_n(M)=0$ for $n<0$.
\end{defn}

By the similar argument as in 
\cite{HA}*{Proposition 7.2.4.11},
we obtain the following proposition. 

\begin{prop}[{cf.~\cite{HA}*{Proposition~7.2.4.11}}]
Let $A$ be a connective Noetherian 
$\mathbb{N}$-graded $\mathbb{E}_{\infty}$-ring.
Then:

\begin{enumerate}[(i)]
\item The full subcategory $\Mod_A^{\afg}(\Z)$
is an idempotent complete small stable
subcategory of $\Mod_A(\Z)$.

\item Every perfect $\Z$-graded $A$-module 
is almost finitely generated.

\item The full subcategory 
$\Mod_A^{\afg}(\Z)_{\ge 0}\subset \Mod_A(\Z)$ 
is closed under the formation of
geometric realizations of simplicial objects.   

\item Let $M$ be an object 
of $\Mod_A^{\afg}(\Z)_{\ge 0}$.
Then $M$ can be obtained as the geometric
realization of a simplicial object $P_{\bullet}$
such that each $P_n$ is a finite direct sum of
elements in $\{A(q)| q\in\mathbb{Z}\}$.
\end{enumerate}
\end{prop}

\begin{defn}\label{ator}
Let $A$ be a connective $\N$-graded $\E$-ring
and let $M$ be 
an almost finitely generated
$\Z$-graded $A$-module.
We say that $M$ is almost torsion
if 
the $\mathbb{Z}$-graded $\pi_0(A)$-module
$\pi_n(M)$ is bounded above
for each $n\in\mathbb{Z}$.
We define an $\infty$-category
\[ \Mod_A^{\ator}(\Z) \]
to be the full subcategory of $\Mod_A^{\afg}(\Z)$
spanned by almost torsion
$\mathbb{Z}$-graded $A$-modules. 
\end{defn}

Next, we construct a truncation
functor
on $\infty$-categories of $\Z$-graded modules over
$\N$-graded $\E$-rings.
Let $i_k: \Z_{\ge k}\hookrightarrow \Z$ be
the inclusion map for $k\in\Z$.
Recall that we have an adjunction
\[ i_{k!}: \Sp(\Z_{\ge k})
   \rightleftarrows
   \Sp(\Z): i_k^*, \] 
of functors between $\infty$-categories,
where the right adjoint $i_k^*$ is the restriction functor
and the left adjoint $i_{k!}$ is 
a left Kan extension along $i_k$.
The $\infty$-categories
$\Sp(\Z)$ and $\Sp(\Z_{\ge k})$ 
are tensored over the
symmetric monoidal $\infty$-category $\Sp(\N)$.
The left adjoint $i_{k!}$ 
is a $\Sp(\N)$-linear functor
and hence the right adjoint $i_{k}^*$ is a lax 
$\Sp(\N)$-linear functor.
By \cite{HA}*{Example~7.3.2.8},
the adjunction $(i_{k!},i_k^*)$ induces
an adjunction
\[ i_{k!}: \Mod_B(\Sp(\Z_{\ge k}))
   \rightleftarrows
   \Mod_B(\Sp(\Z)): i_k^* \] 
for a connective $\N$-graded $\E$-ring $B$.
Note that $\Mod_B(\Sp(\Z))
\simeq \Mod_B(\Z)$.

\begin{defn}
We define a functor
\[ \sigma_{\ge k}: \Mod_B(\Z)\longrightarrow
                  \Mod_B(\Z) \]
to be the composite $\sigma_{\ge k}\simeq i_{k!}\circ i_k^*$.
There is a natural transformation
\[ \epsilon: \sigma_{\ge k}\to {\rm id}, \]
which is a counit of the adjunction $(i_{k!},i_k^*)$.
Since $i_{k!}$ is fully faithful,
$\sigma_{\ge k}$ equipped with $\epsilon$
is a colocalization functor.
\end{defn}

\begin{rem}
For $M\in\Mod_B(\Z)$,
we have $(\sigma_{\ge k}M)_n\simeq M_n$ if $n\ge k$
and $(\sigma_{\ge k}M)_n\simeq 0$ if $n<k$.  
The natural transformation $\epsilon: \sigma_{\ge k}\to {\rm id}$
induces a map
$(\sigma_{\ge k}M)_n\to M_n$
which is an equivalence if $n\ge k$ and
the trivial map if $n<k$.
\end{rem}

\begin{defn}
Let $A$ be a connective $\N$-graded $\E$-ring.
We denote by
\[ \sigma_{\ge k} \Mod_A^{\afg}(\Z) \]
the full subcategory of $\Mod_A^{\afg}(\Z)$
spanned by objects $M$ such that $M_n \simeq 0$ for $n<k$.
We also denote by
\[ \sigma_{\ge k} \Mod_A^{\ator}(\Z) \]
the full subcategory of $\Mod_A^{\ator}(\Z)$
spanned by objects $M$ such that $M_n\simeq 0$ for $n<k$.
\end{defn}

\begin{rem}
If $A$ is Noetherian,
then $\sigma_{\ge k} \Mod_A^{\afg}(\Z)$ 
and $\sigma_{\ge k} \Mod_A^{\ator}(\Z)$
are idempotent complete 
small stable subcategories of $\Mod_A(\Z)$
for any $k\in\mathbb{Z}$.
\end{rem}

Now we study quasi-coherent sheaves
on spectral projective schemes
associated to connective $\N$-graded $\E$-rings.
Let $A$ be a connective Noetherian
$\mathbb{N}$-graded $\mathbb{E}_{\infty}$-ring.
In the following of this subsection we assume that
$A$ satisfies Condition~\ref{condition:N-graded-E-infinity-ring}. 
We let $X=\Proj(A)$ be the spectral projective
scheme associated to $A$.
We denote by $\mathsf{X}$
the underlying ordinary projective scheme $\Proj(\pi_0(A))$
of $X$.
Let $\mathscr{F}$ be a quasi-coherent sheaf of
$\mathcal{O}_{X}$-module on $X$.
For $n\in\mathbb{Z}$,
we have a quasi-coherent sheaf of 
$\mathcal{O}_{\mathsf{X}}$-module 
\[ \pi_n\mathscr{F} \]
on $\mathsf{X}$,
which is the sheafification of the presheaf
obtained by assigning to an open subset $U$ of $\mathsf{X}$
an $\mathcal{O}_{\mathsf{X}}(U)$-module $\pi_n(\mathscr{F}(U))$.

Recall that we have the adjunction
$\widetilde{(-)}:  \Mod_A(\Z) \rightleftarrows
   \QCoh(X): \Gamma_*(X,-)$
of functors between stable $\infty$-categories.

\begin{lemma}\label{lemma:kernel-widetilde-general-case}
For $M\in \Mod_A(\Z)$,
we have $\widetilde{M}\simeq 0$
if and only if
$\pi_n(M)$ is locally bounded above 
for any $n\in\Z$.
\end{lemma}

\begin{proof}
We have 
$\widetilde{M}\simeq 0$ if and only if
$\pi_n\widetilde{M}=0$ for any $n\in\Z$.
We notice that $\pi_n\widetilde{M}$ is equivalent
to the quasi-coherent sheaf on $\mathsf{X}$ associated
to the $\Z$-graded $\pi_0(A)$-module $\pi_n(M)$.
The lemma follows from the fact that 
a quasi-coherent sheaf on $\mathsf{X}$
associated to a $\Z$-graded $\pi_0(A)$-module $N$
is trivial if and only if
$N$ is locally bounded above. 
\end{proof}

\begin{cor}\label{cor:kernel-widetilde-afg-case}
Let $M$ be an almost finitely generated $\Z$-graded $A$-module.
Then $\widetilde{M}\simeq 0$ 
if and only if 
$M$ is almost torsion.
\end{cor}

Now we recall the notion of almost perfect quasi-coherent sheaves
on spectral schemes.

\begin{defn}[{cf.~\cite{HA}*{Proposition~7.2.4.17}}]
Let $R$ be a connective Noetherian $\E$-ring.
We say that an $R$-module $M$ is
almost perfect if 
$\pi_n(M)$ is a finitely generated $\pi_0(R)$-module
for any $n\in\Z$ and if
$\pi_n(M)=0$ for $n\ll 0$.
\end{defn}

According to \cite{SAG}*{Definition~2.8.1.4},
we say that a spectral scheme $Y=(|Y|,\mathcal{O}_Y)$
is locally Noetherian 
if $\mathcal{O}_Y(|V|)$ is a Noetherian $\E$-ring
for any affine open subset $|V|\subset |Y|$.

\begin{defn}[{cf.~\cite{SAG}*{Definition~2.8.4.4}}]  
Let $Y=(|Y|,\mathcal{O}_Y)$ be a locally Noetherian spectral scheme and 
let $\mathscr{F}$ be a quasi-coherent sheaf 
of $\mathcal{O}_Y$-modules on $Y$.
We say that $\mathscr{F}$ is almost perfect if
the restriction $\mathscr{F}|_{U}$
for any affine open subscheme $U\simeq \Spec(R)$
corresponds to an almost perfect $R$-module
under the equivalence
$\QCoh(\Spec(R))\simeq\Mod_R$.
\end{defn}  

Let $\mathscr{F}$ be an almost perfect 
quasi-coherent sheaf of $\mathcal{O}_X$-modules on $X$. 
We consider the quasi-coherent sheaf $\pi_n\mathscr{F}$
of $\mathcal{O}_{\mathsf{X}}$-modules on 
the ordinary projective scheme $\mathsf{X}$. 
First, we shall show that 
$\pi_n\mathscr{F}$ is coherent sheaf on $\mathsf{X}$.

\begin{lemma}\label{lemma:almost-perfect-to-coherent}
If $\mathscr{F}\in \QCoh(X)$
is almost perfect,
then $\pi_n\mathscr{F}$ is a coherent sheaf 
of $\mathcal{O}_{\mathsf{X}}$-modules on 
$\mathsf{X}$ for any $n\in\Z$.
\end{lemma}

\begin{proof}
For a homogeneous element $f$ of $\pi_0(A)$
of positive degree,
we have an affine open subscheme 
$\mathsf{V}\simeq \Spec(\pi_0(A)[f^{-1}]_0)$
of $\mathsf{X}$.
It suffices to show that 
the restriction of $\pi_n\mathscr{F}$
to $\mathsf{V}$ is equivalent to
the quasi-coherent sheaf associated to 
a finitely generated $\pi_0(A)[f^{-1}]_0$-module
for any $f$.

Since $\mathscr{F}$ is almost perfect,
the restriction $\mathscr{F}|_{\mathsf{V}}$ is 
equivalent to the quasi-coherent sheaf associated
to an almost perfect 
$A[f^{-1}]_0$-module $N$.
This implies that $\pi_n\mathscr{F}|_{\mathsf{V}}$
is equivalent to the quasi-coherent sheaf associated to
the $\pi_0(A[f^{-1}]_0)$-module
$\pi_n(N)$.
Since $N$ is almost perfect $A[f^{-1}]_0$-module,
$\pi_n(N)$ is a finitely generated
$\pi_0(A)[f^{-1}]_0$-module.
This completes the proof.
\end{proof}

\begin{lemma}\label{lemma:almost-perfect-bounded-below}
If $\mathscr{F}\in \QCoh(X)$
is almost perfect,
then $\pi_n\mathcal{F}= 0$
for $n\ll 0$.
\end{lemma}

\begin{proof}
We take a set $\{a_i\}_{i=1}^r\subset \pi_0(A_1)$ of generators 
of $\pi_0(A)$ as an $\N$-graded commutative
ring over $\pi_0(A_0)$.
Then we have an affine open covering $\{\mathsf{V}_i\}_{i=1}^r$
of $\mathsf{X}$ given by 
$\mathsf{V}_i\simeq \Spec(\pi_0(A)[a_i^{-1}]_0)$.
By the proof of Lemma~\ref{lemma:almost-perfect-to-coherent},
the restriction of $\pi_n\mathscr{F}$ to $\mathsf{V}_i$
is equivalent to the quasi-coherent sheaf associated to $\pi_n(M_i)$
for some $\Z$-graded $A[a_i^{-1}]_0$-module $M_i$.
Since $\mathscr{F}$ is almost perfect,
$M_i$ is almost perfect and hence
there exists $n_{0}(i)\in\Z$
such that $\pi_n(M_i)= 0$ for $n<n_0(i)$.
Then $\pi_n\mathscr{F}= 0$
for $n<n_0$,
where $n_0=\min\{n_0(i)|\ 1\le i\le r\}$.
\end{proof}

The following proposition
gives a characterization of
the $\Z$-graded $A$-module $\Gamma_*(X,\mathscr{F})$
of global sections for an almost perfect
quasi-coherent sheaf $\mathscr{F}$ of $\mathcal{O}_X$-modules
on a spectral projective scheme $X=\Proj(A)$.

\begin{prop}\label{prop:BKss}
If $\mathscr{F}\in \QCoh(X)$ is almost perfect,
then we have
\begin{itemize}
\item the $\mathbb{Z}$-graded $\pi_0(A)$-module
$\pi_n(\Gamma_*(X,\mathscr{F}))$ is 
strongly quasi-finitely generated
for each $n\in\Z$, and 
\item $\pi_n(\Gamma_*(X,\mathscr{F}))= 0$ for $n\ll 0$.
\end{itemize}
\end{prop}

\begin{proof}
Recall that $\Gamma_*(X,\mathscr{F})$
is defined to be a limit
\[ \lim_{\Delta} \Gamma_*(U_{\bullet},\mathscr{F}) \]
of the cosimplicial object
$\Gamma_*(U_{\bullet},\mathscr{F})$
in $\Mod_A(\Z)$.
We have a Bousfield-Kan spectral sequence
\[ E_2^{s,t}\cong \pi^s\pi_t(\Gamma_*(U_{\bullet},\mathscr{F}))
   \Longrightarrow \pi_{t-s}(\Gamma_*(X,\mathscr{F}))\]
in $\Mod_{\pi_0(A)}^{\heartsuit}(\Z)$ 
abutting to the homotopy group 
of $\Gamma_*(X,\mathscr{F})$
with $d_r: E_r^{s,t}\to E_r^{s+r,t+r-1}$.
Since the cochain complex associated to
the cosimplicial module
$\pi_t\Gamma_*(U_{\bullet},\mathscr{F})$
is isomorphic to the \v{C}ech complex 
$C^{\bullet}(\mathsf{X};\pi_t\mathscr{F}(*))$
of the underlying ordinary projective scheme $\mathsf{X}=\Proj(\pi_0(A))$
with coefficients in the sheaf $\pi_t\mathscr{F}(*)$,
we have an isomorphism
\[ E_2^{s,t}\cong H^s(\mathsf{X};\pi_t\mathscr{F}(*)).\]

Recall that there exists an integer $s_0$
such that $H^s(\mathsf{X};\mathcal{F})=0$
for any quasi-coherent sheaf $\mathcal{F}$ on $\mathsf{X}$
and $s>s_0$.
Hence $E_2^{s,t}=0$ for $s>s_0$.
Furthermore, since $\mathscr{F}$ is almost perfect,
$\pi_n\mathscr{F}(*)= 0$ for $n\ll0$
by Lemma~\ref{lemma:almost-perfect-bounded-below}.
Hence $E_2^{s,t}=0$ for $t\ll 0$.

By Proposition~\ref{prop:cohomology_projective_scheme},
$E_2^{s,t}$ is a strongly quasi-finitely generated
$\mathbb{Z}$-graded $\pi_0(A)$-module
for all $s$ and $t$.
Since $\Mod_{\pi_0(A)}^{\heartsuit\mathrm{sqfg}}(\Z)$
is a Serre subcategory of 
$\Mod_{\pi_0(A)}^{\heartsuit}(\Z)$
by Lemma~\ref{lemma:mod-sqfg-Serre-subcategory},
we see that $E_{\infty}^{s,t}$ is also
strongly quasi-finitely generated for all $s$ and $t$.
Furthermore, $E_2^{s,t}=0$ for $t\ll 0$
implies that 
$E_{\infty}^{s,t}=0$ for $t\ll 0$.

We have a filtration
\[ \pi_n(\Gamma_*(X,\mathscr{F}))=F^{0,n}\supset
   F^{1,n+1}\supset\cdots\supset 
   F^{s_0,n+s_0}\supset 0  \]
such that $F^{s,n+s}/F^{s+1,n+s+1}\cong E_{\infty}^{s,n+s}$.
This implies that 
$\pi_n(\Gamma_*(X,\mathscr{F}))=0$ 
for $n\ll 0$ and 
that
$\pi_n(\Gamma_*(X,\mathscr{F}))$ 
is strongly quasi-finitely generated for all $n$
by using the fact that 
$\Mod_{\pi_0(A)}^{\heartsuit\mathrm{sqfg}}(\Z)$
is a Serre subcategory of 
$\Mod_{\pi_0(A)}^{\heartsuit}(\Z)$.
This completes the proof.
\end{proof}


\subsection{The Serre theorem of almost perfect quasi-coherent
sheaves on projective spectral schemes}
\label{subsection:Serre_theorem_spectral}

In this subsection
we shall prove the Serre theorem 
in spectral algebraic geometry,
which is the main theorem in this paper.

Let $A$ be a connective Noetherian $\N$-graded $\E$-ring.
In this subsection 
we assume that $A$ satisfies
Condition~\ref{condition:N-graded-E-infinity-ring}.
We set $X=\Proj(A)$.
First,
we shall describe the $\infty$-category
$\QCoh(X)$ of quasi-coherent sheaves 
on the projective 
scheme $X$ in terms of $\Z$-graded $A$-modules.

We have the adjunction
$\widetilde{(-)}:\Mod_A(\Z)\rightleftarrows\QCoh(X)
:\Gamma_*(X,-)$
of stable presentable $\infty$-categories,
where the left adjoint $\widetilde{(-)}$ is 
symmetric monoidal 
and the right adjoint $\Gamma_*(-)$
is lax symmetric monoidal and fully faithful.
By Lemma~\ref{lemma:kernel-widetilde-general-case},
we have $\widetilde{M}\simeq 0$
if and only if $\pi_n(M)$ is locally bounded above
for any $n\in \Z$.

\begin{defn}
We say that a $\Z$-graded $A$-module $M$
is locally bounded above in homotopy groups
if the $\Z$-graded $\pi_0(A)$-module $\pi_n(M)$ 
is locally bounded above for each $n\in\Z$.
We define 
\[ \Mod_A^{\mathrm{lbah}}(\Z) \]
to be the full subcategory of $\Mod_A(\Z)$
spanned by those objects that are locally bounded above 
in homotopy groups.
\end{defn}

\begin{rem}
Since $A$ is Noetherian,
$\Mod_A^{\mathrm{lbah}}(\Z)$
is an idempotent complete stable subcategory
of $\Mod_A(\Z)$.
\end{rem}

Let $W$ be the class of all morphisms in $\Mod_A(\Z)$
whose cofiber lies in $\Mod_A^{\mathrm{lbah}}(\Z)$.
We see that a morphism in $\Mod_A(\Z)$
is carried to an equivalence in $\QCoh(X)$
by the functor $\widetilde{(-)}$
if and only if it belongs to $W$.
Hence we obtain the following proposition.

\begin{prop}
The functor
$\widetilde{(-)}: \Mod_A(\Z)\to\QCoh(X)$
induces an equivalence
\[ \Mod_A(\Z)/\Mod_A^{\mathrm{lbah}}(\Z)\stackrel{\simeq}
   \longrightarrow \QCoh(X) \]
of stable symmetric monoidal presentable $\infty$-categories,
where the left hand side is 
the localization with respect to the class $W$.
\end{prop}

Next, we consider the problem
to describe the $\infty$-category $\QCoh(X)^{\aperf}$
of almost perfect quasi-coherent sheaves on $X$
in terms of $\Z$-graded $A$-modules.
For this purpose,
we recall Verdier quotients
of small stable $\infty$-categories
discussed in \cite{BGT}*{\S5} and \cite{NS}*{\S{I.3}}.

For stable $\infty$-categories $\mathcal{A}$ and $\mathcal{B}$,
we denote by $\Fun^{\mathrm{ex}}(\mathcal{A},\mathcal{B})$
the $\infty$-category of exact functors from $\mathcal{A}$
to $\mathcal{B}$.
Let $\mathcal{C}$ be a small stable $\infty$-category 
and let $\mathcal{D}$ be a stable subcategory 
of $\mathcal{C}$.
By \cite{NS}*{Theorem~I.3.3},
the Verdier quotient 
\[  \mathcal{C}/\mathcal{D} \]
is a small stable $\infty$-category
equipped with an exact functor
$\mathcal{C}\to \mathcal{C}/\mathcal{D}$
which has the following universal property:
For any stable $\infty$-category $\mathcal{E}$,
the composition with $\mathcal{C}\to\mathcal{C}/\mathcal{D}$
induces an equivalence between 
$\Fun^{\mathrm{ex}}(\mathcal{C}/\mathcal{D},\mathcal{E})$
and the full subcategory of
$\Fun^{\mathrm{ex}}(\mathcal{C},\mathcal{E})$
spanned by those functors which carry
all objects of $\mathcal{D}$ to $0$.

We also have a multiplicative property
of the construction of Verdier quotients.
For stable symmetric monoidal $\infty$-categories
$\mathcal{A}$ and $\mathcal{B}$,
we denote by $\Fun^{\mathrm{ex}}_{\otimes}(\mathcal{A},\mathcal{B})$
the $\infty$-category of 
symmetric monoidal exact
functors from $\mathcal{A}$ to $\mathcal{B}$.
We suppose that $\mathcal{C}$
is a small stable symmetric monoidal $\infty$-category
in which the tensor product functor
$\otimes: \mathcal{C}\times\mathcal{C}\to\mathcal{C}$
is exact separately in each variable.
We say that a full subcategory $\mathcal{D}$
is a $\otimes$-ideal if $X\otimes Y\in\mathcal{D}$
for any $X\in\mathcal{C}$ and $Y\in\mathcal{D}$.
In this case,
by \cite{NS}*{Theorem~I.3.6},
the Verdier quotient $\mathcal{C}/\mathcal{D}$
acquires a unique symmetric monoidal structure 
such that $\mathcal{C}\to\mathcal{C}/\mathcal{D}$ is 
a symmetric monoidal exact functor,
and it has the following universal property:
The composition with $\mathcal{C}\to\mathcal{C}/\mathcal{D}$
induces an equivalence between
$\Fun^{\mathrm{ex}}_{\otimes}(\mathcal{C}/\mathcal{D},\mathcal{E})$
and
the full subcategory of 
$\Fun^{\mathrm{ex}}_{\otimes}(\mathcal{C},\mathcal{E})$
spanned by those functors which carry
all objects of $\mathcal{D}$ to $0$. 

Now, we compare the $\infty$-category
$\QCoh(X)^{\aperf}$ of almost perfect quasi-coherent
sheaves on $X$ with
the $\infty$-category
$\Mod_A^{\afg}(\Z)$
of almost finitely generated $\Z$-graded $A$-modules.

The $\infty$-category $\QCoh(X)^{\aperf}$ is
symmetric monoidal
since it is a full subcategory of $\QCoh(X)$
which contains the unit object $\mathcal{O}_X$ and is closed
under tensor products.
The $\infty$-category
$\Mod_A^{\afg}(\Z)$ is also
a symmetric monoidal $\infty$-category.

\begin{lemma}\label{lemma:mod-A-afg-to-QCoh-X-aperf}
The functor 
$\widetilde{(-)}:\Mod_A(\Z)\to\QCoh(X)$
induces
a symmetric monoidal exact functor
\[ \widetilde{(-)}: \Mod_A^{\afg}(\Z)\longrightarrow\QCoh(X)^{\aperf} \]
by restriction.
\end{lemma}

\begin{proof}
It suffices to show that the quasi-coherent sheaf
$\widetilde{M}$ is almost perfect if
a $\Z$-graded $A$-module $M$ is almost finitely generated.
Let $f\in\pi_0(A)$ be homogeneous of positive degree
and let $U=\Spec(A[f^{-1}]_0)$ be
an affine open subscheme of $X$.
The restriction $\widetilde{M}|_U$
corresponds to the $A[f^{-1}]_0$-module $M[f^{-1}]_0$.
We see that $M[f^{-1}]_0$ is almost perfect since
$M$ is almost finitely generated.
This completes the proof.
\end{proof}


\begin{lemma}\label{lemma:symmetric-monoidal-widehat}
The functor
$\widetilde{(-)}: \Mod_A^{\afg}(\Z)\to \QCoh(X)^{\aperf}(X)$
induces 
a symmetric monoidal exact functor
\[ \widehat{(-)}: \Mod_A^{\afg}(\Z)/\Mod_A^{\ator}(\Z)
   \longrightarrow \QCoh(X)^{\aperf}.\] 
\end{lemma}

\begin{proof}
The tensor product functor of 
$\Mod_A^{\afg}(\Z)$ 
is exact separately in each variable.
For $M\in\Mod_A^{\afg}(\Z)$,
by Corollary~\ref{cor:kernel-widetilde-afg-case},
we have $\widetilde{M}\simeq 0$ if and only if
$M\in\Mod_A^{\ator}(\Z)$.
In particular, 
we see that $\Mod_A^{\ator}(\Z)$
is a $\otimes$-ideal of $\Mod_A^{\afg}(\Z)$.
The lemma follows from 
the universal property of Verdier quotients.
\end{proof}

In order to prove that $\widehat{(-)}$
is an equivalence,
we shall compare $\QCoh(X)^{\aperf}$
with the truncation $\sigma_{\ge k}\Mod_A^{\afg}(\Z)$
of the $\infty$-category of almost finitely generated
$\Z$-graded $A$-modules for any $k\in\Z$.

For each $k\in {\mathbb{Z}}$,
we define a functor
\[ \Gamma_{\ge k}(X,-): \QCoh(X)\longrightarrow
                  \Mod_A(\Z) \]
to be the composite 
$\Gamma_*(X,-): \QCoh(X)\to \Mod_A(\Z)$
with $\sigma_{\ge k}:  \Mod_A(\Z) \to \Mod_A(\Z)$. 
By Proposition~\ref{prop:BKss},
we see that 
$\Gamma_{\ge k}(X,\mathscr{F})$
is an object of $\sigma_{\ge k} \Mod_A^{\afg}(\Z)$
if $\mathscr{F}\in \QCoh(X)$ is almost perfect.

\if
If $M\in \Mod_A(\Z)$ is almost finitely generated,
then the quasi-coherent sheaf 
$\widetilde{M}\in \QCoh(X)$ associated to $M$
is almost perfect.
\fi
On the other hand,
we have a functor 
\[ \widetilde{(-)}: \sigma_{\ge k} \Mod_A^{\afg}(\Z)
   \longrightarrow \QCoh(X)^{\aperf}\]
by Lemma~\ref{lemma:mod-A-afg-to-QCoh-X-aperf}.

\begin{prop}\label{prop:equivalence-aperf-mod-A-ge-k}
For each $k\in \Z$,
we have an adjunction
\[ \widetilde{(-)}: \sigma_{\ge k} \Mod_A^{\afg}(\Z)
   \rightleftarrows
   \QCoh(X)^{\aperf}: \Gamma_{\ge k}(X,-)\]
of $\infty$-categories.
The right adjoint
$\Gamma_{\ge k}(X,-)$ is fully faithful.
The left adjoint $\widetilde{(-)}$
induces an equivalence
\[ \sigma_{\ge k} \Mod_A^{\afg}(\Z)/
   \sigma_{\ge k} \Mod_A^{\ator}(\Z)
   \stackrel{\simeq}{\longrightarrow}
   \QCoh(X)^{\aperf} \]
of small stable $\infty$-categories.
\end{prop}

\begin{proof}
Since $\sigma_{\ge k}$ is a colocalization functor,
we have an equivalence of mapping spaces
\[ \Map_{\Mod_A(\Z)}(M,\Gamma_{\ge k}(X,\mathscr{F}))
   \stackrel{\simeq}{\longrightarrow}
   \Map_{\Mod_A(\Z)}(M,\Gamma_*(X,\mathscr{F}))\]
if $M\in\Mod_A(\Z)$ satisfies
$\epsilon: \sigma_{\ge k}M\stackrel{\simeq}{\to} M$.
This implies that
\[ \Map_{\sigma_{\ge k} \Mod_A^{\afg}(\Z)}
    (M,\Gamma_{\ge k}(X,\mathscr{F}))
   \stackrel{\simeq}{\longrightarrow}
   \Map_{\QCoh(X)^{\aperf}}(\widetilde{M},\mathscr{F}),\]
and hence the pair $(\widetilde{(-)},\Gamma_{\ge k}(X,-))$
is an adjunction of $\infty$-categories.
Since the counit $\widetilde{(-)}\circ \Gamma_{\ge k}(X,-)\to 
{\rm id}$ of the adjunction is an 
equivalence,
the right adjoint $\Gamma_{\ge k}(X,-)$ is fully faithful.
By Corollary~\ref{cor:kernel-widetilde-afg-case},
for $M\in\sigma_{\ge k} \Mod_A^{\afg}(\Z)$,
we have $\widetilde{M}\simeq 0$
if and only if 
$M\in\sigma_{\ge k} \Mod_A^{\ator}(\Z)$.
This implies an equivalence of $\infty$-categories
between the Verdier quotient 
$\sigma_{\ge k} \Mod_A^{\afg}(\Z) /
\sigma_{\ge k} \Mod_A^{\ator}(\Z) $ 
and $\QCoh(X)^{\aperf}$.
\end{proof}

We have a commutative diagram
\[ \begin{array}{ccc}
    \vdots & & \vdots\\
    \downarrow & & \downarrow \\[2mm]
    \sigma_{\ge -k} \Mod_A^{\afg}(\Z) /
   \sigma_{\ge -k} \Mod_A^{\ator}(\Z)
   &\stackrel{\simeq}{\longrightarrow}&
   \QCoh(X)^{\aperf}\\[2mm]
   \mbox{$\rm \scriptstyle{\simeq}$}\bigg\downarrow 
   \phantom{\mbox{$\rm \scriptstyle{\simeq}$}}
   & & 
   \phantom{\mbox{$\rm \scriptstyle{=}$}}\bigg\downarrow
   \mbox{$\rm \scriptstyle{=}$}\\[2mm]
    \sigma_{\ge -(k+1)} \Mod_A^{\afg}(\Z) /
   \sigma_{\ge -(k+1)} \Mod_A^{\ator}(\Z)
   &\stackrel{\simeq}{\longrightarrow}&
   \QCoh(X)^{\aperf}\\[2mm]
    \downarrow & & \downarrow \\[2mm]
    \vdots & & \vdots\\
   \end{array}\] 
of small stable $\infty$-categories.
By Proposition~\ref{prop:equivalence-aperf-mod-A-ge-k}, 
the horizontal arrows are equivalences 
induced by the functor $\widetilde{(-)}$.
The left vertical arrows are induced by the inclusion functors. 

\begin{thm}\label{n3}
The functor $\widetilde{(-)}: \Mod_A(\Z)\to \QCoh(X)$
induces an equivalence
\[ \Mod_A^{\afg}(\Z)/ \Mod_A^{\ator}(\Z)
   \stackrel{\simeq}{\longrightarrow}
   \QCoh(X)^{\aperf} \]
of small stable symmetric monoidal  $\infty$-categories.
\end{thm}


\begin{proof}
Since we have 
the symmetric monoidal exact functor
\[ \widehat{(-)}: \Mod_A^{\afg}(\Z)/\Mod_A^{\ator}(\Z)
   \longrightarrow \QCoh(X)^{\aperf} \]
by Lemma~\ref{lemma:symmetric-monoidal-widehat},
it suffices to show that
it is an equivalence of small stable $\infty$-categories.

We fix any integer $k$.
We recall that the functor $\widetilde{(-)}:
\Mod_A^{\afg}(\Z)\to \QCoh(X)^{\aperf}$
is the composite of 
$i_k^*: \Mod_A^{\afg}(\Z)\to
\sigma_{\ge k} \Mod_A^{\afg}(\Z)$
with
$\widetilde{(-)}:\sigma_{\ge k} \Mod_A^{\afg}(\Z)
\to \QCoh(X)^{\aperf}$
and that the functor 
$\widetilde{(-)}:\sigma_{\ge k} \Mod_A^{\afg}(\Z)
\to \QCoh(X)^{\aperf}$
induces an equivalence of small stable $\infty$-categories 
between 
$\sigma_{\ge k} \Mod_A^{\afg}(\Z)/
\sigma_{\ge k} \Mod_A^{\ator}(\Z)$
and $\QCoh(X)^{\aperf}$.
Thus, it suffices to show that $i_k^*$ induces
an equivalence 
between
$\Mod_A^{\afg}(\Z)/ \Mod_A^{\ator}(\Z)$
and
$\sigma_{\ge k} \Mod_A^{\afg}(\Z)/
\sigma_{\ge k} \Mod_A^{\ator}(\Z) $.

We have an adjunction
\[ i_{k!}: \sigma_{\ge k}\Mod_A^{\afg}(\Z)\rightleftarrows
      \Mod_A^{\afg}(\Z): i_k^*\]
of $\infty$-categories,
where $i_{k!}$ is the inclusion functor.
The functor $i_{k!}$ induces 
a functor   
\[ \widehat{i_{k!}}: \sigma_{\ge k}\Mod_A^{\afg}(\Z)/
   \sigma_{\ge k} \Mod_A^{\ator}(\Z)\longrightarrow
   \Mod_A^{\afg}(\Z)/ \Mod_A^{\ator}(\Z) . \]
For $M\in \Mod_A^{\afg}(\Z)$,
we have $i_k^*M\in \sigma_{\ge k}\Mod_A^{\ator}(\Z)$
if and only if $M\in \Mod_A^{\ator}(\Z)$.
Hence the functor $i_k^*$
induces a functor
\[   \widehat{i_k^*}: \Mod_A^{\afg}(\Z)/ 
     \Mod_A^{\ator}(\Z)
     \longrightarrow 
     \sigma_{\ge k} \Mod_A^{\afg}(\Z)/
     \sigma_{\ge k} \Mod_A^{\ator}(\Z). \]
It is clear that $\widehat{i_k^*}\circ \widehat{i_{k!}}
\simeq{\rm id}$.
For any $M\in \Mod_A^{\afg}(\Z)$,
the cofiber of the counit map
$(i_{k!}\circ i_k^*)(M)\to M$ belongs to 
$\Mod_A^{\ator}(\Z)$.
This implies that 
$\widehat{i_{k!}}\circ\widehat{i_k^*}\simeq \mathrm{id}$.
Hence $\widehat{i_k^*}$
gives an equivalence of stable $\infty$-categories.
\end{proof}

\begin{rem}\label{bbafgtor}
We say that $M\in\Mod_A(\Z)$ is 
bounded below 
if there exists $n_0\in\Z$
such that 
$M_n\simeq 0$ for $n<n_0$. 
We define an $\infty$-category
\[ \Mod_A^{\afg}(\Z)_{-} \]
to be the full subcategory of $\Mod_A^{\afg}(\Z)$
spanned by bounded below almost finitely generated 
$\Z$-graded $A$-modules.
We also define an $\infty$-category
\[ \Mod_A^{\ator}(\Z)_{-} \]
to be the full subcategory of $\Mod_A^{\ator}(\Z)$
spanned by bounded below almost torsion 
$\Z$-graded $A$-modules.
In the same way as 
the proof of Theorem~\ref{n3},
we
can show that the functor $\widetilde{(-)}$
induces an equivalence 
\[ \Mod_A^{\afg}(\Z)_{-}/ \Mod_A^{\ator}(\Z)_{-}
   \stackrel{\simeq}{\longrightarrow}
   \QCoh(X)^{\aperf} \]
of small stable symmetric monoidal $\infty$-categories.
\end{rem}

\end{document}